\documentclass[12pt]{amsart}
\theoremstyle{plain}

\newtheorem{theorem}{Theorem}
\newtheorem{definition}{Definition}
\newtheorem{corollary}[theorem]{Corollary}

\newtheorem{example}[theorem]{Example}
\newtheorem{lemma}[theorem]{Lemma}
\newtheorem{proposition}[theorem]{Proposition}
\newtheorem{remark}[theorem]{Remark}
\newtheorem{claim}[theorem]{Claim}

\def\endex{{\hfill\rule{7pt}{1pt} }}

\def\id{1\!\!1}
\def\BOX{\raisebox{1.2mm}{\hskip .5mm $\fbox{\hphantom{\hglue .01mm}}$}
         \hskip .5mm}
\def\C{{\mathcal C}}
\def\ot{\otimes} 
\def\op{\oplus}

\def\boldk{{\bf k}} \def\bfk{\boldk} \def\otexp#1#2{#1^{\otimes #2}}
\def\D{{\sf D}}
\def\C{{\sf C}}
\def\P{{\mathcal P}}
\def\calP{{\mathcal P}}
\def\calL{{\mathcal L}}
\def\mathcalr{{\mathcal R}}
\def\Ass{{{\it Ass\/}}}
\def\Comm{{{\it Comm\/}}}
\def\Lie{{{{\it Lie\/}}}}
\def\Comm{{{{\it Comm\/}}}}
\def\calM{{\mathcal M}}
\def\coll{{\mbox{\tt Coll\/}}}
\def\Coll#1{\{#1(n)\}_{n\geq 1}}
\def\Set#1{\{#1\}}
\def\oper{{\mbox{\tt Oper\/}}}
\def\k{{\bf k}}
\def\kVect{\mbox{${\tt gr}$-${\tt Vect}_{\bf k}$}}
\def\Sets{{\tt Sets}}
\def\AbGrp{{\mbox{{\tt Ab}-{\tt Grp}}}}
\def\S{{{\mathcal S}}} \def\bold{\bf}
\def\I{{\mathcal I}}
\def\F{{\mathcal F}}
\def\pmod{{\mbox{$\P$-{\tt Mod\/}}}}
\def\freep#1{{\P\langle#1\rangle}}
\def\prez#1#2{\langle #1;#2 \rangle}
\def\susp{\uparrow\!}
\def\desusp{\downarrow\!}
\def\Ker#1{{\mbox{\rm Ker}(#1)}}
\def\Coker#1{{\mbox{\rm Coker}(#1)}}
\def\free#1#2{{#1\langle#2\rangle}}
\def\calR{{\mathcal R}}
\def\Span{{\mbox{\rm Span}}}

\def\H{{\mathcal H}} \def\O{{\mathcal O}}
\def\V{{\mathcal V}} 

\def\T{{\mathcal T}}
\def\sign#1{{(-1)^{{#1}}}}

\def\G{{\sf G}}

\def\jac#1#2#3{[#1[#2#3]]+[#2[#3#1]]+[#3[#1#2]]}

\def\bb{\bullet}
\def\cc{\circ}
\def\K{{\mathcal K}}
\def\setdots#1#2{\{#1,\ldots,#2\}}
\def\D{{\mathcal D}}
\def\ix#1#2#3{1#1(2#2(3#34))}
\def\ii#1#2#3{(1#12)#2(3#34)}
\def\iii#1#2#3{((1#12)#23)#34}
\def\iv#1#2#3{(1#1(2#23))#34}
\def\vv#1#2#3{1#1((2#23)#34)}
\def\BB{\!\bullet\!}
\def\CC{\!\circ\!}
\def\Hom{{\mbox{\rm Hom}}}
\def\HOM#1#2{{\rm Hom}_W(#1,#2)}
\def\End#1{{\mbox{\rm End}(#1)}}

\def\homa-a#1#2{{\mbox{\rm Hom}_{A}(#1,#2)}}
\def\homa-v#1#2{{\mbox{\rm Hom}_{V}(#1,#2)}}

  \def\ker{{\mbox{\rm Ker}}}

\def\OO{\!\otimes\!}
\def\rada#1#2{{#1,\ldots,#2}}
\def\Rada#1#2#3{#1_{#2},\dots,#1_{#3}}
\def\orada#1#2{{#1 \otimes \ldots \otimes #2}}
\def\cases#1#2#3#4{
                  \left\{
                         \begin{array}{ll}
                           #1,\ &\mbox{#2}
                           \\
                           #3,\ &\mbox{#4}
                          \end{array}
                   \right.
}

\def\sgn{{\mbox{\rm sgn}}}
\def\bk{{\bf k}}
\def\L{{\mathcal L}}

\def\Vn{{V^{\otimes n}}}
\def\eps{{\epsilon}}
\def\lra{\longrightarrow} 
\def\sta{\stackrel}
\def\epic{{\to \hskip -3mm \to}}
\def\De{\Delta}
\def\verylra{
\unitlength=1em
\thinlines
      \begin{picture}(3,.5)(0,-.3)
      \vector(1,0){3}  
      \end{picture}
}

\def\pentagon#1#2#3#4#5{
\begin{center}
\begin{picture}(4,4)(-2,-0.5)
\thicklines

\put(-1,0){\vector(1,0){1.9}}
\put(1,0){\vector(1,2){.95}}
\put(0,3){\vector(2,-1){1.9}}
\put(-2,2){\vector(2,1){1.9}}
\put(-2,2){\vector(1,-2){.95}}

\put(1,0){\makebox(0,0){$\bullet$}}
\put(2,2){\makebox(0,0){$\bullet$}}
\put(0,3){\makebox(0,0){$\bullet$}}
\put(-2,2){\makebox(0,0){$\bullet$}}
\put(-1,0){\makebox(0,0){$\bullet$}}

\put(1,-.5){\makebox(0,0)[l]{{$#1$}}}
\put(2.5,2){\makebox(0,0)[l]{{$#2$}}}
\put(0,3.5){\makebox(0,0){{$#3$}}}
\put(-2.5,2){\makebox(0,0)[r]{{$#4$}}}
\put(-1,-.5){\makebox(0,0)[r]{{$#5$}}}

\end{picture}
\end{center}}

\def\pentagonline#1#2#3#4#5{
\begin{center}
\begin{picture}(4,4)(-2,-0.5)
\thicklines

\put(-1,0){\line(1,0){1.9}}
\put(1,0){\line(1,2){.95}}
\put(0,3){\line(2,-1){1.9}}
\put(-2,2){\line(2,1){1.9}}
\put(-2,2){\line(1,-2){.95}}

\put(1,0){\makebox(0,0){$\bullet$}}
\put(2,2){\makebox(0,0){$\bullet$}}
\put(0,3){\makebox(0,0){$\bullet$}}
\put(-2,2){\makebox(0,0){$\bullet$}}
\put(-1,0){\makebox(0,0){$\bullet$}}

\put(1,-.5){\makebox(0,0)[l]{{$#1$}}}
\put(2.5,2){\makebox(0,0)[l]{{$#2$}}}
\put(0,3.5){\makebox(0,0){{$#3$}}}
\put(-2.5,2){\makebox(0,0)[r]{{$#4$}}}
\put(-1,-.5){\makebox(0,0)[r]{{$#5$}}}

\end{picture}
\end{center}}

\newcommand{\square}[8]{
\setlength{\unitlength}{1.01cm}
\begin{picture}(5,3.6)
\thicklines

\put(0,3){\makebox(0,0){$#1$}}
\put(5,3){\makebox(0,0){$#2$}}
\put(0,0){\makebox(0,0){$#3$}}
\put(5,0){\makebox(0,0){$#4$}}

\put(-.5,1.5){\makebox(0,0)[r]{$#6$}}
\put(5.5,1.5){\makebox(0,0)[l]{$#7$}}
\put(2.5,-0.5){\makebox(0,0)[t]{$#8$}}
\put(2.5,3.5){\makebox(0,0)[b]{$#5$}}

\put(1,0){\vector(1,0){3}}
\put(1,3){\vector(1,0){3}}
\put(0,2.5){\vector(0,-1){2}}
\put(5,2.5){\vector(0,-1){2}}
\end{picture}
}

\def\ssquare#1#2#3#4{
\unitlength=1.00mm
\thicklines
\begin{picture}(46.00,35.00)
\put(40.00,10.00){\line(-1,0){30.00}}
\put(10.00,10.00){\line(0,1){20.00}}
\put(10.00,30.00){\line(1,0){30.00}}
\put(40.00,30.00){\line(0,-1){20.00}}
\put(10.00,30.00){\makebox(0,0)[cc]{$\bullet$}}
\put(40.00,30.00){\makebox(0,0)[cc]{$\bullet$}}
\put(10.00,10.00){\makebox(0,0)[cc]{$\bullet$}}
\put(40.00,10.00){\makebox(0,0)[cc]{$\bullet$}}
\put(5.00,35.00){\makebox(0,0)[cc]{$#1$}}
\put(5.00,5.00){\makebox(0,0)[cc]{$#3$}}
\put(45.00,5.00){\makebox(0,0)[cc]{$#4$}}
\put(46.00,35.00){\makebox(0,0)[cc]{$#2$}}
\end{picture}
}

\def\hex#1#2#3#4#5#6{
\unitlength=1mm
\thicklines
\begin{picture}(55.00,35.00)
\put(50.00,10.00){\line(-1,0){40.00}}
\put(10.00,10.00){\line(0,1){20.00}}
\put(10.00,30.00){\line(1,0){40.00}}
\put(50.00,30.00){\line(0,-1){20.00}}
\put(10.00,30.00){\makebox(0,0)[cc]{$\bullet$}}
\put(30.00,30.00){\makebox(0,0)[cc]{$\bullet$}}
\put(50.00,30.00){\makebox(0,0)[cc]{$\bullet$}}
\put(10.00,10.00){\makebox(0,0)[cc]{$\bullet$}}
\put(30.00,10.00){\makebox(0,0)[cc]{$\bullet$}}
\put(50.00,10.00){\makebox(0,0)[cc]{$\bullet$}}
\put(5.00,35.00){\makebox(0,0)[cc]{$#1$\hskip2mm}}
\put(30.00,35.00){\makebox(0,0)[cc]{$#2$}}
\put(55.00,35.00){\makebox(0,0)[cc]{\hskip2mm$#3$}}
\put(5.00,5.00){\makebox(0,0)[cc]{$#4$\hskip2mm}}
\put(30.00,5.00){\makebox(0,0)[cc]{$#5$}}
\put(55.00,5.00){\makebox(0,0)[cc]{\hskip2mm$#6$}}
\end{picture}
}

\begin{document}
\bibliographystyle{plain}
\pagestyle{plain}

\title{COHERENCE CONSTRAINTS FOR OPERADS, CATEGORIES AND ALGEBRAS}

\author[M. Markl]{Martin Markl}
\author[S. Shnider]{Steve Shnider}

\thanks{M.~Markl was partially 
        supported by the Academic Exchange Program and by the
        grant GA AV \v CR 1019507. 
        During part of the time this research was done, S.~Shnider
        was partially supported by The Israel Science Foundation 
        founded by The Israel Academy of Sciences and Humanities -- 
        Center of Excellence Program. 
}

\catcode`\@=11
\address{M.~M.: Mathematical Institute of the Academy\\ 
         \v Zitn\'a 25\\
         115 67 Praha 1\\ \hfill\break The~Czech~Republic}
\email{markl@math.cas.cz}
\address{S.~S.: Department of Mathematics\\
         University of Bar Ilan, Israel}
\email{shnider@macs.biu.ac.il}
\catcode`\@=13



\begin{abstract}
Coherence phenomena appear in two different situations. 
In the context of category theory the term `coherence constraints'
refers to a set of diagrams whose commutativity implies the
commutativity of a larger class of diagrams. In the context of algebra
coherence constrains are a minimal set of generators for the second
syzygy, that is, a set of equations which generate the full set of
identities among the defining relations of an algebraic theory.

A typical example of the first type
is Mac Lane's coherence theorem for monoidal 
categories~\cite[Theorem~3.1]{maclane:RiceUniv.Studies63}, an example
of the second type is the result of~\cite{drinfeld:Alg.iAnaliz89} 
saying that
pentagon identity for the `associator' $\Phi$ of a quasi-Hopf algebra 
implies the validity of a set of identities with higher instances of $\Phi$. 

We show that both types of coherence are
governed by a homological invariant of the operad
for the underlying algebraic structure. We call this invariant 
the (space of) coherence
constraints.  In many cases these constraints can be
explicitly described, thus giving rise to various coherence results, both
classical and new.
\end{abstract}

\maketitle

\section{Introduction}
\label{1}
\begin{figure}[ht]
\def\pentagon1#1#2#3#4#5{
\begin{center}
\begin{picture}(4,4)(-2,-0.5)
\thicklines

\put(-1,0){\vector(1,0){1.9}}
\put(1,0){\vector(1,2){.95}}
\put(0,3){\vector(2,-1){1.9}}
\put(-2,2){\vector(2,1){1.9}}
\put(-2,2){\vector(1,-2){.95}}

\put(1,0){\makebox(0,0){$\bullet$}}
\put(2,2){\makebox(0,0){$\bullet$}}
\put(0,3){\makebox(0,0){$\bullet$}}
\put(-2,2){\makebox(0,0){$\bullet$}}
\put(-1,0){\makebox(0,0){$\bullet$}}

\put(1,-.5){\makebox(0,0)[l]{{$#1$}}}
\put(2.5,2){\makebox(0,0)[l]{{$#2$}}}
\put(0,3.5){\makebox(0,0){{$#3$}}}
\put(-2.5,2){\makebox(0,0)[r]{{$#4$}}}
\put(-1,-.5){\makebox(0,0)[r]{{$#5$}}}

\put(-1.6,2.9){\makebox(0,0){{$a\circ (\id^2 \times \BOX)$}}}
\put(1.6,2.9){\makebox(0,0){{$a\circ (\BOX \times \id^2)$}}}
\put(0,.3){\makebox(0,0){{%
      $a\! \circ \!(\id\!\! \times%
\!\! \BOX \!\!\times\!\! \id)$}}}
\put(-2,1){\makebox(0,0)[r]{{$\BOX\circ (\id \times a)$}}}
\put(2,1){\makebox(0,0)[l]{{$\BOX\circ(a \times \id)$}}}

\end{picture}
\end{center}}

\setlength{\unitlength}{1.2cm}
\pentagon1{\BOX(\BOX\times \id)(\id \times \BOX\times \id)}%
{\BOX(\BOX\times \id)(\BOX\times \id^2)}%
{\BOX(\BOX \times \BOX)}%
{\BOX(\id \times \BOX)(\id^2\times \BOX)}%
{\BOX(\id\times \BOX)(\id \times \BOX\times \id)}
\caption{\label{pentagon}
The Pentagon}
\end{figure}
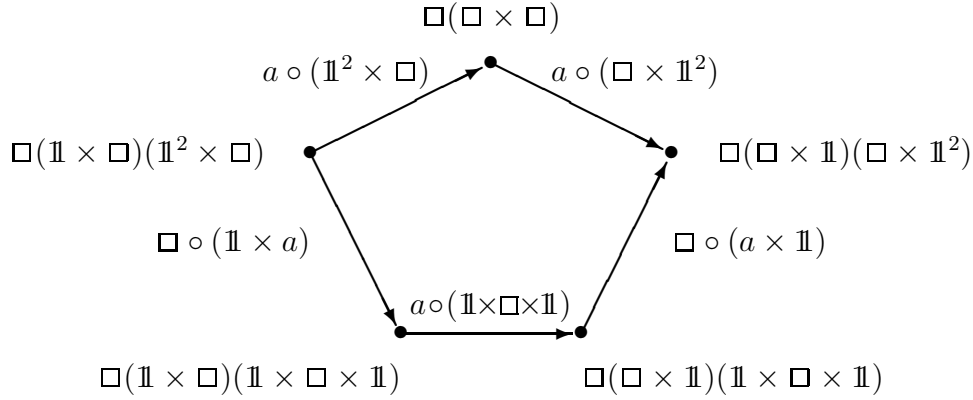

We remind the reader of some definitions and results
of~\cite{maclane:RiceUniv.Studies63}.
A {\em category with a multiplication\/} is
a category ${\mathcal C}$ together with
a covariant bifunctor $\BOX : {\mathcal C} \times
{\mathcal C} \to {\mathcal C}$. An {\em associativity
isomorphism\/} for
$({\mathcal C},\BOX)$ is then a natural transformation
\begin{equation}
\label{ass-iso}
a : \BOX(\id \times \BOX) \longrightarrow \BOX(\BOX \times \id)
\end{equation}
($\id$ denotes the identity functor)
such that each $a(A,B,C) : A\BOX (B\BOX C)\to (A\BOX B)\BOX C$
has a two-sided inverse in ${\mathcal C}$,
for $A,B,C\in {\mathcal C}$; here we denote, as usual,
$\BOX(\id \times \BOX)(A,B,C)$ by $A\BOX (B\BOX C)$,~etc.
Having such an associativity isomorphism, we can
consider diagrams whose vertices are iterates of $\BOX$ and edges
expansions of instances of $a$. The category ${\mathcal C}$ is called
{\em coherent\/} if all these diagrams commute.
The easiest of these diagrams is the
pentagon (see Figure~\ref{pentagon}).
There is no a~priori reason for the commutativity of any of these diagrams,
but the celebrated Mac~Lane's coherence 
theorem~\cite[Theorem~3.1]{maclane:RiceUniv.Studies63} says that the
commutativity of one diagram, the pentagon, 
implies the commutativity of all these
diagrams. 

Consider the case of a $\bfk$-vector
space $U$ which is a module over a unital associative, not necessary
coassociative bialgebra $V= (V,\cdot,\Delta,1)$. 
We adopt the Drinfel'd convention and consider the associativity
$a^{-1}$ represented by the action of an invertible element
$\Phi = \sum_i
\mbox{$\Phi_{1,i} \!\ot\! \Phi_{2,i} \!\ot\! \Phi_{3,i}$}\in V^{\ot
3}$. Let $* : U\ot
U \to U$ be a bilinear product which is `$\Phi$-associative,'
\begin{equation}
\label{Ubu}
\Phi(a*(b*c)) = ((a*b)*c),
\mbox {  for $a,b,c\in U$,}
\end{equation}
where $\Phi(a*(b*c))$ is an abbreviation for 
$\sum (\Phi_1a*(\Phi_2b*\Phi_3c))$. 
Assuming $(U,*)$ is a
$V$-module algebra, \cite[Chapter~10]{Majid:book}, that is
\[
v\cdot (a*b) = v_{(1)}\cdot a *v_{(2)}\cdot b,
\mbox { for $v\in V$, $a,b,c \in U$,}
\]
(we will follow usual conventions and delete the 
summation sign and summation indices use the Sweedler 
abbreviated notation with $v_{(1)}\otimes v_{(2)}$ 
standing for $\Delta(v)=\sum_i v_{(1)i}\otimes v_{(2)i}$). 
One derives easily from~(\ref{Ubu}) that
\begin{equation}
\label{Majka_a_Bobecek}
(((a*b)*c)*d) =P(\Phi)(((a*b)*c)*d),
\end{equation}
where
$P(\Phi):=(\Phi \ot 1)^{-1}(\id \ot \De \ot \id)
(\Phi)^{-1}(1\ot \Phi)^{-1}(\id \ot \De)
(\Phi)(\De \ot \id)(\Phi).$

If we do not assume that $P(\Phi) = 1$, 
then~(\ref{Majka_a_Bobecek}) is a new relation 
imposed on the space of bracketed 4-fold products other than the 
association relation to other bracketings. The
condition $P(\Phi) = 1$ is, of course, the same as
\begin{equation}
\label{pentagon_equation}
(1\ot\Phi)(\id \ot \Delta \ot\id)(\Phi)
(\Phi\ot 1)=(\id^2 \ot \Delta)(\Phi)(\Delta \ot \id^2)(\Phi), 
\end{equation}
the famous pentagon condition on $\Phi$ introduced by 
Drinfel'd~\cite{drinfeld:Alg.iAnaliz89}, although not from this point
of view.
It is called so because its five factors correspond to the five sides
of the pentagon for the natural monoidal structure on the category of
modules over the algebra $(V,\cdot)$.

We may consider
$*$-products of order $5$ and higher and look for 
similar equations in $\Phi$. Because of Mac Lane's theorem, 
all these equations follow
from the pentagon condition~(\ref{pentagon_equation}). 

We presented two situations where the coherence appears -- 
one in category theory, where it was formulated in
terms of commutative diagrams, and another in algebra, where it was
expressed in the language 
of algebraic equations. Both examples above were related to a
certain associativity -- 
of the transformation $\BOX$ in~(\ref{ass-iso}) and  of
the product $*$ in~(\ref{Ubu}). Also the description of the `coherence
constraint' pointed in both cases to the shape of the pentagon.
In this case the algebraic result was derived from the category
theory. We also 
want to show how to derive categorial results from algebra.

For an operad $\P$, we introduce, in
Definition~\ref{61}, the space
$\C_{\P}$ of coherence constraints of $\P$. It is a certain
homological invariant of $\P$ that can be informally described
as a second syzygy, spanned by `the relations among the
defining relations' where the `defining relations' generate
the ideal defining  an operad as a quotient of a
free operad. These coherence constraints can be read off from the
bigraded model of $\P$ (Corollary~\ref{csa}) 
and can be easily described for so
called Koszul operads (Theorem~\ref{main}).

We show that both types of coherence boil down to
coherence constraints of the governing operad. For the categorial
coherence, 
it is Theorem~\ref{21} which says, roughly speaking, that
the commutativity of diagrams corresponding to $\C_{\P}$ implies the
commutativity of all diagrams. 

The algebraic situation is more subtle. The 
coherence of a `quantization' is given in terms of
linear equations and it intuitively means that any solution of a
linear equation in the
original system deforms to a solution of the quantized one.
 
This is formalized by introducing the `$V$-relative' version $\P_V$ of
the operad $\P$,
the coherence then means that $\P_V$ is a `flat extension' of $\P$
(Definition~\ref{32}),
that is, the defining relations of the quantized structure have the
same rank as the original structure.
The characterization of this kind of coherence
in terms of coherence constraints of $\P$ is given in
Theorem~\ref{Tatouch}.

We illustrate our methods by giving a new 
proof of Mac Lane's coherence ($\P = \Ass$, the operad for associative
algebras), which is just
two lines once we know that $\Ass$ is Koszul (see Example~\ref{25}).

Example~\ref{[]} is related to the linear logic
and Example~\ref{digebres} presents a categorial version of Loday's
bigebras. 
  
On the quantum side we discuss Drinfel'd's quasi-Hopf algebras
(Example~\ref{cervena_propiska}, 
$\P = \Ass$), generalized Lie algebras analogous to the construction
given by
Gurevich (Example~\ref{modra_propiska}, 
$\P = \Lie$) and a certain form of
strictly homotopy associative algebras (Example~\ref{nozik}).

As an useful tool, we introduce in Section~\ref{4}  
a series of  bipartite graphs, 
encoding the coherence relations. 
In many cases these relations can be described by
a simple graph which we call  the {Tel-A-graph\/}
(Tel-A= Tel~Aviv, the place where the discovery was made). While, for
$\P=\Ass$ this Tel-A-graph is the pentagon, for $\P=\Lie$ it is the
Peterson graph (Figure~\ref{peterson}), 
and for the operad governing algebras of Example~\ref{nozik} 
a kind of the M\"obius strip!

\noindent 
{\it Acknowledgment.\/}
The authors would like to express their gratitude to Jim Stasheff for
reading the manuscript and many useful remarks and comments.

\section{Operads and coherence constraints}
\label{2}

The notion of an {\em operad\/} and of an {\em algebra\/} over an operad is
classical and well known (see~\cite{may:1972}, or more recent
sources~\cite{ginzburg-kapranov:DMJ94,%
getzler-jones:preprint,markl:zebrulka}).
We thus recall only briefly the definitions and notation. We will
need also some results on homology and presentations of operads; this
part of the paper relies on~\cite{markl:zebrulka}.

Operads make sense in any strict symmetric monoidal
category $\mathcal M = (\calM,\BOX)$. 
The most important example for the purposes of this paper is
the category \kVect\ of graded ${\bf k}$-vector spaces,
where ${\bf k}$ is a field of characteristic zero, with monoidal structure
given by the standard graded tensor product over ${\bf k}$. 
In Sections 4 and 5 we consider
also the category \Sets\  of sets, with monoidal structure given by the 
cartesian product, and the category \AbGrp\ of abelian groups with tensor
product over ${\mathbb Z}$.
 
More precisely, an operad is a sequence $\P = \{\P(n); n\geq 1\}$ of
objects of $\calM$
such that:
\begin{itemize}
\item[(i)]
Each $\P(n)$ is equipped with a (right) action of the
symmetric group $\Sigma_n$ on $n$ elements, $n\geq 1$.
\item[(ii)]
For any $m_1,\ldots,m_l \geq 1$ we have the {\em composition maps}
$$
\gamma=\gamma_{m_1,\ldots,m_l}:\P(l)\BOX\P(m_1)\BOX\cdots\BOX
\P(m_l) \longrightarrow\P(m_1+\cdots+m_l).
$$
\end{itemize}
These data have to satisfy the usual axioms including the existence
of a unit $1\in \P(1)$, for which we refer
to~\cite{may:1972}. We sometimes write $\mu(\nu_1,\cdots,\nu_l)$,
$\mu(\nu_1 \BOX \cdots \BOX \nu_l)$ or $\gamma(\mu;\nu_1,\cdots,\nu_l)$
instead of $\gamma(\mu\BOX\nu_1\BOX\cdots\BOX\nu_l)$ (notice that in
all three above examples $(\kVect,\ot)$, $(\AbGrp,\ot)$ and
$(\Sets,\times)$ of the monoidal category $(\calM,\BOX)$ 
the `product of elements'
$\mu\BOX\nu_1\BOX\cdots\BOX\nu_l$ makes sense).

A {\em collection\/} is a sequence $E = \{E(n);\ n\geq 2\}$
of elements of $\calM$ such that each $E(n)$ is equipped with an
action of the symmetric group $\Sigma_n$.
The obvious forgetful functor ${F}_{\rm or} :\oper \to
\coll$ from the category of
operads to the category of collections
has a left adjoint $\F
:\coll \to \oper$ and we call
$\F(E)$ the {\em free operad\/} on the collection $E$.

\begin{remark}{\rm
The notions above, as well as all the results which follow, have
obvious {\em non-$\Sigma$\/} (also called {\em nonsymmetric\/})
analogs which we obtain by forgetting
everything related to the symmetric group action. We thus have
non-$\Sigma$ operads, non-$\Sigma$ collections, etc. The reason for
considering these objects is that the non-$\Sigma$ versions are much
simpler and there are many examples which live
in a non-$\Sigma$ world, the most prominent being the
associative algebra case. We will move freely between $\Sigma$
and non-$\Sigma$ worlds clarifying when necessary
the context in which we are working.\endex
}\end{remark}

In the rest of
this section, an operad is an operad in the 
category $\kVect$ of graded ${\bf k}$-vector spaces.

A {\em module\/} over an operad $\P$ is an abelian group
object in the slice category $\oper/\P$. The axioms were explicitly
given for modules over a so-called {\em pseudo-operad}
in~\cite{markl:zebrulka}, in the standard case the axioms are quite
analogous. Namely, a module over $\P$ is a collection
$M = \{M(n);\ n\geq 1\}$ together with a map
\begin{eqnarray*}
&
m:\displaystyle\bigoplus_{1\leq i\leq l}\{
\P(l)\!\ot\!\P(m_1)\!\ot\!\cdots\!\ot\!
M(m_i) \!\ot\!
\P(m_l)\} \!\oplus\!\{ M(l)\!\ot\!\P(m_1)\!\ot\!\cdots\!\ot\!
\P(m_l) \}
\longrightarrow&
\\
&\hskip2cm\longrightarrow
M(m_1+\cdots+m_l)&
\end{eqnarray*}
given for any $m_1,\ldots,m_l \geq 1$. This map is supposed to
satisfy obvious axioms given by the linearization of the axioms of
operads. 
Just as for the operadic composition map, we sometimes
write $a(\rada{b_1}{b_l})$, $a(b_1 \otimes \cdots \otimes b_l)$ 
instead of $m(a\otimes b_1\otimes \cdots
\otimes b_l)$.

We will give some examples of $\P$-modules which we will need in the
sequel. The operad $\P$ itself is a $\P$-module. If
$\alpha :\P \to \S$ is an operad map, then $\alpha$ induces a
$\P$-module structure on $\S$. Finally, if $\I \subset \P$ is an
ideal in $\P$ (see~\cite{markl:zebrulka}), then $\I$ is naturally a
$\P$-module.

The forgetful functor $F_{\rm or} :\pmod \to \coll$ from the category of
$\P$-modules to the category of collections has a left adjoint
$\freep-:\coll \to \pmod$, and we call the
$\P$-module
$\freep E$ the {\em free $\P$-module\/} on the collection $E$.

{}For the time being, we suppose that our operads always have $\P(1)=
{\bf k}$.
Let $\P^+\subset \P$ be the ideal defined by $\P^+(1):=0$ and
$\P^+(n):=\P(n)$ for $n\geq 2$. 

For a $\P$-module $(M,m)$ we define the {\em decomposables\/} 
of $M$ to be the collection $D(M) = D_\P(M)$
generated by elements either of the form $m(r;p_1,\ldots,p_l)$, 
where $r\in M$, $p_1,\ldots,p_l \in P$ and at 
least one of $p_1,\ldots,p_l\in \P$ belongs to
$\P^+$, 
or of the form
$m(p;p_1,\ldots,r,\ldots,p_l)$, where again $r\in M$, $p,p_1,\ldots,p_l\in
\P$ and at least one of $p,p_1,\ldots,p_l\in \P$ belongs to
$\P^+$. Define
the {\em indecomposables\/} of $M$ as the collection $Q(M) =
Q_\P(M):= M/D_\P(M)$. 
  
Each operad can be represented as $\P = \F(E)/(R)$, where $E$ and $R$
are collections and $(R)$ is the operadic ideal generated by $R$; we
write $\P = \prez ER$. Because $\P(1)= {\bf k}$, we can always
suppose that the presentation is minimal~\cite{markl:zebrulka}.
This means, by definition, that $E \cong Q_\P(\P^+)$ and that the
collection $R$ is isomorphic to the indecomposables of the kernel of
the canonical map $\F(E)\to \P$, $R \cong Q_{\F(E)}\{\Ker{\F(E)\to \P}\}$.

Following Quillen's paradigm, we 
consider the higher derived functors of the functor of indecomposables, 
see~\cite{quillen:LNM43}.

Let $\P = \prez ER$. Let $J:= \free{\F(E)}{R}$ be the free
$\F(E)$-module on $R$ and let $\pi: J \to (R)$
be the obvious natural epimorphism
of $\F(E)$-modules. For $x\in \free{\F(E)}{R}(l)$,
$y\in \free{\F(E)}{R}$ and
$a_1,\ldots,a_l\in \F(E)$ the element
\[
o' := x(a_1,\ldots,a_{s-1},\pi(y),a_{s+1},\ldots,a_l) -
\pi(x)(a_1,\ldots,a_{s-1},y,a_{s+1},\ldots,a_l) \in J
\]
belongs to $\Ker{\pi}$, for any $1\leq s\leq l$. Similarly, for
$b\in \F(E)(l)$, $a_1,\ldots,a_l\in \F(E)$ and
$x,y\in \free{\F(E)}{R}$, the element
\begin{eqnarray*}
\lefteqn{
o'' :=
b(a_1,\ldots,a_{s-1},\pi(x),a_{s+1},
\ldots,a_{t-1},y,a_{t+1},\ldots,a_l) -\hskip 2cm}
\\
&&\hskip 2cm -b(a_1,\ldots,a_{s-1},x,a_{s+1},
\ldots,a_{t-1},\pi(y),a_{t+1},\ldots,a_l) \in J
\end{eqnarray*}
belongs to $\Ker \pi$, for any $1\leq s<t\leq l$. In the spirit of
the definition of the cotangent cohomology we call the $\F(E)$-module
generated by elements of the above two types the {\em module
of obvious
relations\/} and denote it by $\O=\O_\P$. To understand better the
meaning of this module we recommend 
looking at the
definition of the module $U_0$ on page~44
of~\cite{lichtenbaum-schlessinger:TAMS67}
in the classical commutative
algebra situation or to the definition in 2.2 of~\cite{markl:ws93}.

\begin{definition}
\label{61}
The collection of {\em coherence relations} of the operad $\P$ is the
$\F(E)$-module $\D = \D_\P := \Ker{\pi:J \to (R)}/\O$. The
collection of {\em coherence constraints} of the operad $\P$ is the collection of
indecomposables of the $\F(E)$-module $\D$,
$\C = \C_\P:= Q_{\F(E)}(\D)$.
\end{definition}

\begin{example}{\rm\
\label{hedges}
Let $\xi$ be an independent symbol and let $E$ be the non-$\Sigma$
collection defined by
\[
E(n): =
\cases{\Span(\xi)}{for $n=2$, and}{0}{otherwise.}
\]
Let $r:= \xi(1,\xi)- \xi(\xi,1) \in \F(E)(3)$ and let
$R$ be the
(non-$\Sigma$) collection generated by $r$, i.e.
\[
R(n):=
\cases{\Span(r)}{for $n=3$, and}0{otherwise.}
\]
Associative algebras are
then the algebras over the non-$\Sigma$ operad $\Ass := \prez ER$.

There are five rooted planar binary trees having four leaves so
$\dim_{\boldk}(\F(E)(4))= 5$. Moreover, there are five rooted
planar trees with one bivalent and one trivalent vertex and four leaves
so $\dim_{\boldk}(\free{\F(E)}R)(4) =5$.
 Choose as basis elements  for $\F(E)(4)$ the trees  represented
in the standard way by the following bracketings
\begin{eqnarray*}
&{\bf a}:= ((12)3)4,\ {\bf b}:= (12)(34),\ {\bf c}:= 1(2(34)),&
\\
&{\bf d}:= 1((23)4),\ {\bf e}:=(1(23))4 \in \F(E)(4).&
\end{eqnarray*}
Choose also the basis
\[
{\bf 1}:= r(\xi,1,1),\ {\bf 2}:= r(1,1,\xi),\ {\bf 3}:= \xi(1,r),\
{\bf 4}:= r(1,\xi,1),\ {\bf 5}:= \xi(r,1) \in \free{\F(E)}R,
\]
for $\F(E)\langle R \rangle(4)$.
The map $\pi: \F(E)\langle R \rangle(4) \to (R)(4)$
has the following matrix description:
\[
\begin{tabular}{r|rrrrr}
               &${\bf a}$&${\bf b}$&${\bf c}$&${\bf d}%
                 $&${\bf e}$\\  \hline
$\pi({\bf 1}) $&$-1 $&$1 $&$0 $&$0 $&$0 $\\
$\pi({\bf 2}) $&$ 0 $&$-1$&$1 $&$ 0$&$0 $\\
$\pi({\bf 3}) $&$ 0 $&$0 $&$1$&$-1 $&$0 $\\
$\pi({\bf 4}) $&$ 0 $&$ 0$&$ 0$&$1 $&$-1$\\
$\pi({\bf 5}) $&$-1 $&$ 0$&$0 $&$ 0$&$+1$\\
\end{tabular}
\]
One sees immediately that $\dim(\Ker \pi)(4)=1$, and that the
kernel is spanned by the element $p\in J(4) = \free{\F(E)}{R}(4)$ defined by
\begin{equation}
\label{assrelate}
p:= \xi(r,1) -r(\xi,1,1)+ r(1,\xi,1)- r(1,1,\xi)+\xi(1,r).
\end{equation}
The five terms in $p$ correspond to the edges of the pentagon
considered as the Stasheff associahedron $K_4$ with five vertices
labeled by the five binary trees with four leaves or,
equivalently, by the five possible bracketings of four elements.
We may prove, by a step-by-step repeating the arguments of the proof
of~\cite[Theorem~3.1]{maclane:RiceUniv.Studies63}, that
$p$ generates the collection of coherence constraints $\C$, as was
in fact done in the last section of~\cite{markl:ws93}, but we
derive this statement using a more sophisticated approach,
which we now describe.
\endex}\end{example}
\noindent

\section{Coherence and the homology of operads}
\label{3a}

In this section we work with operads in the monoidal category 
\kVect\ of graded vector spaces. Let us recall some notions and results
of~\cite[Section~3]{markl:zebrulka}.

Let $\S$ be an operad. By a {\em differential\/} on $\S$ we mean a
degree $-1$ map $d: \S \to \S$ of collections having the expected
Leibniz
property with respect to the operadic composition
and satisfying $d^2=0$. A differential on the free
operad $\F(E)$ is uniquely determined by its restriction to the space
of generators $E$.

Suppose we have a collection $Z$ which decomposes as $Z=Z^0\op
Z^1\op\cdots$ (meaning, of course, that for each $n\geq 2$ we have a
$\Sigma_n$-invariant decomposition $Z(n) = Z^0(n)\op
Z^1(n)\op\cdots$ of the graded vector space $Z(n) = \bigoplus Z_j(n)$
into the direct sum of graded vector spaces $Z^k(n) = \bigoplus
Z^k_j(n)$, $k \geq 0$).
This induces on $\F(Z)$ a grading, $\F(Z) =
\bigoplus_{k\geq 0}\F(Z)^k$. We call this grading the
{\em TJ-grading\/}
(from Tate-Jozefiak) here; the reason will became obvious below.

In the situation above, the free operad $\F(Z)$ is bigraded,
$\F(Z) = \bigoplus \F(Z)^k_j$, where  $k$ refers to the
TJ-grading introduced above and
$j$ indicates the `inner' grading given by the grading of
$Z = \bigoplus Z_j$.

Suppose that $d$ is an (inner) degree $-1$ differential on $\F
(Z)$ such
that
\begin{equation}
\label{3}
d(Z^k) \subset \F(Z)^{k-1}, \mbox{ for all $k\geq 1$}
\end{equation}
(meaning, of course, that $d(Z^k)(n) \subset \F(Z)^{k-1}(n)$ for
all
$n\geq 2$), i.e.~that $d$ is homogeneous degree $-1$ with respect to
the
TJ-grading. Then the homology operad $\H(\F(Z),d)$ is
bigraded,
\[
\H(\F(Z),d) = \bigoplus \H^k_j( \F (Z),d),
\]
the upper grading being induced by the TJ-grading and the lower
one by
the inner grading. In~\cite{markl:zebrulka} the first author
proved the following theorem.

\begin{theorem}
\label{TJ}
Let ${\mathcal P}$ be an operad (\hglue1mm$\Sigma$ or non-$\Sigma$,
with trivial differential). Then
there
exists a collection $Z=Z^0\op
Z^1\op\cdots$, a differential $d$ on $ \F (Z)$
satisfying~(\ref{3}) and
a map $\rho : (\F(Z),d)\to ({\mathcal P},0)$ of differential operads
such that the following conditions
are satisfied:
\begin{enumerate}
\item[(i)]
$d$ is minimal in the sense that $d(Z)$ consists of decomposable
elements of the operad $\F(Z)$,
\item[(ii)]
$\rho|_{Z^{\geq 1}}=0$ and $\rho$ induces an isomorphism
$\H^0(\rho):\H^0(\F(Z),d) \cong {\mathcal P}$, and
\item[(iii)]
$\H^{\geq 1}(\F(Z),d)=0$.
\end{enumerate}
\end{theorem}
We call the object $\rho :(\F(Z),d)\to({\mathcal P},0)$ 
 the {\em bigraded model\/} of the operad ${\mathcal P}$.
This object is an analog
of the bigraded model of a commutative
graded algebra constructed in~\cite[Section~3]{halperin-stasheff}.

Recall that the {\em suspension\/} of a graded vector space
$V = \bigoplus V_i$ is the graded vector space $\susp V$ defined by
$(\susp V)_i = V_{i-1}$.
The following proposition is obvious from the construction of the
bigraded model described in~\cite{markl:zebrulka}.

\begin{proposition}
Let $\P$ be an operad and let $\P = \prez ER$ be a minimal
presentation. Let $(\F(Z),d)$ be a bigraded model of $\P$. Then there
are the following isomorphisms of collections:
\begin{equation}
\label{dymka}
Z^0 \cong E,\ Z^1 \cong \susp R,\ Z^2 \cong \susp^2 \C
\end{equation}
\end{proposition}

The most important for our purpose is the third equation
of~(\ref{dymka}). We formulate the following corollary.

\begin{corollary}
\label{csa}
The collection $\C_\P$ of coherence constraints of an operad $\P$ is
isomorphic to the double desuspension $\desusp^2 Z^2$ of the
collection
of TJ-degree two indecomposables of the bigraded model of $\P$.
\end{corollary}

We say that an operad $\P$ is {\em quadratic\/} if it has a
presentation $\P = \prez ER$ such that $E(n)=0$ for $n\not=2$ and
$R(n)=0$ for $n\not=3$. Each quadratic operad has its {\em quadratic
dual\/} $\P^!$~\cite{ginzburg-kapranov:DMJ94},
which is another operad constructed very
explicitly from the presentation of $\P$. V.~Ginzburg and
M.M.~Kapranov
introduced in~\cite{ginzburg-kapranov:DMJ94} an extremely important
notion of the Koszulness of an operad. It is a certain homological
property analogous to the Koszulness of commutative algebras; operads
sharing this property are called Koszul operads.
The following proposition appeared in \cite{markl:zebrulka}
as Proposition 2.6.

\begin{proposition}
\label{nuzky}
Let $\P$ be a quadratic operad and suppose that $\P$ is Koszul.
Denote by $\P^!$ its quadratic dual. Then there is, for each $n\geq 1$,
the following isomorphism of $\Sigma_n$-modules:
\[
Z^i(n)\cong
\left\{
\begin{array}{ll}
0,\ & n\not=i+2
\\
\sgn \otimes\susp^{(n-2)} \P^!(n),\ & n=i+2,
\end{array}
\right.
\]
where $\sgn$ denotes the one-dimensional signum representation.
For non-$\Sigma$ operads, $Z^i(n)$ is given by
the same formula but
without the signum representation.
\end{proposition}

\noindent 
Combining Proposition~\ref{nuzky} with Corollary~\ref{csa} we get
the following proposition.

\begin{theorem}
\label{main}
For a Koszul quadratic operad $\P$, there a $\Sigma_4$-equivariant
isomorphism 
\[
\C_\P = \C_P(4)\cong \sgn \otimes \P^!(4).
\]
In particular, the only coherence constraints are in degree~4.
\end{theorem}

Let us formulate a proposition counting the dimension of the
space of coherence constraints.
We need the generating function $g_\P(x)$ of an operad $\P$,
which is the formal power series
\begin{equation}
\label{Carolina_Inn}
g_\P(x) :=
\sum_{n\geq 1}\frac{\dim(\P(n))}{n!}x^n.
\end{equation}
As it follows
from~\cite{ginzburg-kapranov:DMJ94}, if $\P$ is Koszul, then
\[
g_\P(-g_{\P^!}(-x))=x
\]
The same formula holds also for a
non-$\Sigma$ $\P$ if we drop the $n!$ (= the order of $\Sigma_n$) 
from~(\ref{Carolina_Inn}). The formula above enables one to express,
for a Koszul operad $\P$, the
dimension $\dim(\P^!(4))$ via $\dim(\P(2))$, $\dim(\P(3))$ and
$\dim(\P(4))$:

\begin{proposition}
\label{za}
Suppose that $\P$ is a non-$\Sigma$ Koszul operad. Then
\[
\dim(\C_\P) = \dim(\P(4)) +5\dim(\P(2))[\dim(\P(2))^2-\dim(\P(3))].
\]
If $\P$ is symmetric, then
\[
\dim(\C_\P) = \dim(\P(4)) +5\dim(\P(2))[3\dim(\P(2))^2-2\dim(\P(3))].
\]
\end{proposition}

\begin{example}{\rm\
\label{Lie}
Let $\zeta$ be an independent variable and $E$ the symmetric
collection defined by
\[
E(n):=
\cases{\Span(\zeta)}{for $n=2$, and}0{otherwise,}
\]
with the sign representation of $\Sigma_2$ on $E(2)$.
Let $B_n$ be the free nonassociative anticommutative algebra on
the set
$\{1,\ldots,n\}$ and let $B'_n$ denote the subset of $B_n$ spanned by
monomials in which each element of $\{1,\ldots,n\}$ appears
exactly once. Then $\F(E)(n)\cong B'_n$ for any $n\geq 1$, where
$\F(E)$ is now  the free {\em symmetric\/} operad on the
collection $E$.

Let $\iota := \jac123 \in \F(E)(3)$ and $R$ be the symmetric
collection
generated by $\iota$. Then $\Lie := \prez ER$ is the operad
governing  Lie algebras.
This operad is quadratic 
Koszul~\cite[page~229]{ginzburg-kapranov:DMJ94} and
$\Lie^! =
\Comm$, therefore, by Theorem~\ref{main}, $\C = \Comm(4)$, which
is the 
one dimensional trivial representation of $\Sigma_4$.
Thus there is only one coherence constraint,
as in the associative algebra case. Let us describe this constraint
explicitly.  

A description of the map $\pi :\F(E)\langle R
\rangle(4) \to (R)(4)$ is given by the matrix of Figure~\ref{matrix}.   
\def\aa{{\bf a}}
\def\bb{{\bf b}}
\def\cc{{\bf c}}
\def\dd{{\bf d}}
\def\ee{{\bf e}}
\def\ff{{\bf f}}
\def\gg{{\bf g}}
\def\hh{{\bf h}}
\def\ii{{\bf i}}
\def\jj{{\bf j}}
\def\kk{{\bf k}}
\def\ll{{\bf l}}
\def\mm{{\bf m}}
\def\nn{{\bf n}}
\def\oo{{\bf o}}
\def\I{{\bf 1}}
\def\II{{\bf 2}}
\def\III{{\bf 3}}
\def\IV{{\bf 4}}
\def\V{{\bf 5}}
\def\VI{{\bf 6}}
\def\VII{{\bf 7}}
\def\VIII{{\bf 8}}
\def\IX{{\bf 9}}
\def\X{{\bf 10}}
\begin{figure}[ht]
\[
\begin{tabular}{r|rrrrr|rrrrrrrrrr}
${}$&$\aa$&$\bb$&$\cc$&$\dd$&%
\multicolumn{1}{r}{$\ee$}&$\ff$&$\gg$&$\hh$&$\ii$&$\jj$&$\kk$&%
$\ll$&$\mm$&$\nn$&$\oo$
\\  \hline
$\pi(\I)$&$+1$&$-1$&$ 0$&$ 0$&$ 0$&$ 0$&$ 0$&$ 0$&$ 0$&$
0$&$0$&$+1$&%
$ 0$&$ 0$&$ 0$
\\
$\pi(\II)$&$ 0$&$+1$&$-1$&$ 0$&$ 0$&$ 0$&$+1$&$ 0$&$ 0$&$ 0$&%
$ 0$&$ 0$&$ 0$&$ 0$&$ 0$
\\
$\pi(\III)$&$ 0$&$ 0$&$-1$&$+1$&$0$&$+1$&$ 0$&$ 0$&$ 0$&$ 0$&%
$ 0$&$ 0$&$ 0$&$ 0$&$ 0$
\\
$\pi(\IV)$&$ 0$&$ 0$&$ 0$&$+1$&$-1$&$ 0$&$ 0$&$ 0$&$ 0$&$ 0$&$%
 0$&$ 0$&$ 0$&$ 0$&$+1$
\\
$\pi(\V)$&$-1$&$ 0$&$ 0$&$ 0$&$+1$&$ 0$&$ 0$&$ 0$&$ 0$&$ 0$&$%
0$&$ 0$&$+1$&$ 0$&$ 0$
\\ \cline{2-6}
$\pi(\VI)$&$ 0$&$ 0$&$ 0$&$ 0$&\multicolumn{1}{r}{$ 0$}&%
$ 0$&$ 0$&$ 0$&$-1$&$ 0$&$%
0$&$ 0$&$+1$&$-1$&$ 0$
\\
$\pi(\VII)$&$ 0$&$ 0$&$ 0$&$ 0$&%
\multicolumn{1}{r}{$ 0$}&$-1$&$ 0$&$ 0$&$ 0$&$+1$&$%
0$&$ 0$&$ 0$&$+1$&$ 0$
\\
$\pi(\VIII)$&$ 0$&$ 0$&$ 0$&$ 0$&%
\multicolumn{1}{r}{$ 0$}&$ 0$&$-1$&$+1$&$-1$&$ 0$&%
$ 0$&$ 0$&$ 0$&$ 0$&$ 0$
\\
$\pi(\IX)$&$ 0$&$ 0$&$ 0$&$ 0$&%
\multicolumn{1}{r}{$ 0$}&$ 0$&$ 0$&$ 0$&$ 0$&$-1$&%
$+1$&$+1$&$ 0$&$ 0$&$ 0$
\\
$\pi(\X)$&$ 0$&$ 0$&$ 0$&$ 0$&%
\multicolumn{1}{r}{$ 0$}&$ 0$&$ 0$&$-1$&$ 0$&$ 0$&%
$+1$&$ 0$&$ 0$&$ 0$&$-1$
\\
\end{tabular}
\]
\noindent
The basis elements are chosen as:
\def\arraystretch{1.3}
\[
\begin{tabular}{rrrrr}
$\aa:=[[[12]3]4]$&$\bb:=[[12][34]]$&$\cc:=
[1[2[34]]]$&$\dd:= [1[[23]4]]$&$\ee := [[1[23]]4]$
\\
$\ff:= [1[3[24]]]$&$\gg:=[2[1[34]]]$&$\hh
:=[2[3[14]]]$&$\ii:=[2[4[13]]]$&$\jj:=[3[1[24]]]$
\\
$\kk:=[3[2[14]]]$&$\ll:=[3[[12]4]]$&$\mm
:=[4[2[13]]]$&$\nn:= [[13][24]]$&$\oo:= [[14][23]]$
\\
\end{tabular}
\]
\noindent
and:
\[
\begin{tabular}{lll}
${\I}: = -\iota(\zeta,1,1)$&
${\II}:= -\iota(1,1,\zeta)$&
${\III}: = -\zeta(1,\iota)$
\\
${\IV}:=- \iota(1,1,\zeta)\cdot T_{1342}$&
${\V}:= - \zeta(1,\iota)\cdot T_{2341}$&
${\VI}:= -\iota(\zeta,1,1)\cdot T_{1324}$
\\
${\VII}:= -\iota(1,1,\zeta)\cdot T_{1324}$&
${\VIII}:= -\zeta(1,\iota)\cdot T_{2134}$&
${\IX}:= -\zeta(1,\iota)\cdot T_{2314}$
\\
$\hskip -2.3mm {\X}:= -\iota(\zeta,1,1)\cdot T_{1342}$&
{}&{}
\end{tabular}
\]
\caption{
\label{matrix}
The matrix of the map $\pi : {\mathcal F}(E)\langle R \rangle(4) \to
(R)(4)$. The symbol $T_{i_1i_2i_3i_4}$ denotes 
the permutation which
sends $(1234)$ to $(i_1i_2i_3i_4)$.
The upper left 5$\times$5-submatrix coincides, up to sign,
with the corresponding matrix for the associative algebra operad.
}
\end{figure}
Let $\ell \in \free{\F(E)}{R}(4)$ be the
element defined as
\begin{eqnarray*}
\ell &:=&
\iota(\zeta,1,1)
+ \iota(1,1,\zeta)
-\zeta(1,\iota)
+\iota(1,1,\zeta)\cdot T_{1342}
+\zeta(1,\iota)\cdot T_{2341}
\\
&&-\iota(\zeta,1,1)\cdot T_{1324}
-\iota(1,1,\zeta)\cdot T_{1324}
+\zeta(1,\iota)\cdot T_{2134}
\\
&&-\zeta(1,\iota)\cdot T_{2314}
+\iota(\zeta,1,1)\cdot T_{1342},
\end{eqnarray*}
or, in the notation of Figure~\ref{matrix},
\[
\ell = -{\bf 1} -{\bf 2} +{\bf 3} 
-{\bf 4} -{\bf 5} +{\bf 6} +{\bf 7} -{\bf 8} +{\bf 9} -{\bf 10}.
\]
It is easy to check that $\ell \in
\Ker\pi$.
By Theorem~\ref{main}, $\ell$ is the only coherence constraint
for Lie algebras.
\endex}\end{example}

\section{A graphic description of coherence relations}
\label{4}

A nice feature of the associative
operad  is that the coherence relations can be represented
by commutative diagrams, the closed edge-paths in the 
one skeleton of the Stasheff associahedra~\cite{stasheff:TAMS63}. 
In this section we show how  to describe the coherence relations 
for any quadratic operad by bipartite graphs, and in certain cases 
by ordinary graphs.

Associated to any presentation of a quadratic operad $\P = \prez ER$ there
is a series of bipartite graphs,
$\T= \{\T(n)\}_{n\geq 4}$. Recall that $J$ denoted
the free $\F(E)$-module on $R$, $\pi : J\to (R)$ was the canonical
epimorphism and $\D := \Ker\pi/\O$, see Section~2.
Let $n\geq 4$. Fix a basis $(v_1,\ldots,v_b)$
(resp.~$(r_1,\ldots,r_s)$) of $\F(E)(n)$ (resp.~of
$J(n)=\free{\F(E)}R(n)$).
Then define a bipartite graph with vertices partitioned into the  two sets
${\mathcal V}(n):= \{v_1,\ldots,v_b\}$ and ${\mathcal V}'(n) := \{r_1,\ldots,r_s\}$.
The vertices $v_i\in {\mathcal V}(n)$ and $r_j\in {\mathcal V}'(n)$ are joined by 
an edge if and only if $v_i$ occurs with a non-zero coefficient in $\pi(r_j)$.

There is a simple graph encoding the same data whenever 
for one of the sets ${\mathcal V}(n)$ or  ${\mathcal V}'(n)$ there are
precisely two edges incident to each of the vertices in that set.
For example, if ${\mathcal V}'(n)$ has this property, we create a new graph by
concatenating the two edges incident to  a vertex of ${\mathcal V}'(n)$, 
delete the vertex and create a single edge connecting 
two vertices in ${\mathcal V}(n)$ and give the new edge the same label as 
the deleted vertex. The new graph has edges labeled by the elements of 
${\mathcal V}'(n)$ and vertices labeled by the elements of ${\mathcal V}(n)$. 
In this case each edge is labeled by a basis element in $\free{\F(E)}R(n)$
relating the  two elements of the basis of $\F(E)(n)$
which label the endpoints.
In this case we say that the defining relations are
{\em graphlike\/} and call the corresponding graph $\G(\T(n))$
the {\em Tel-A-graph} (the idea arose from discussions in Tel Aviv).

Note that the bipartite graph $\T(n)$ does not  encode all the data
from the presentation of $\P=\prez ER$ of an operad
since it only shows which  terms $v_i$ appear  
in $\pi(r_j)$ with nonzero coefficient, 
not the actual values of the coefficients.
On the other hand if we can choose the bases $\{v_i\}$ for all the $\F(E)(n)$
such that the coefficients are all $\pm 1$ then all 
the data of the presentation
can be encoded into the graph $\T(n)$ by orienting the edges.
The edge connecting 
the vertex $v_i\in {\mathcal V}(n)$ and the relation $r_j\in {\mathcal V}'(n)$ 
will be oriented in the direction of $v_i$ if the coefficient 
of $v_i$ in the relation $r_j$ is $+1$ and away from $v_i$ if 
the coefficient in $-1$. 
If $\T(n)$ is graphlike and oriented, that does not
guarantee that the Tel-A-graph has a consistent 
orientation. 
What is required is that
there exist an orientation such that each vertex with two incident edges in the
bipartite graph has one incoming edge and one outgoing edge. 
In this case whenever
such a vertex is deleted, the new edge created has a natural orientation.

We will see that the existence of an oriented Tel-A-graph
formalizes the property of an operad being defined by axioms expressed by 
commutative diagrams. 
{}For instance, 
the operad {\it Ass} has an oriented
Tel-A-graph $\G(\T(n))$ given by the one skeleton of the 
Stasheff associahedron $K_n$, with orientation as shown, for $n=4$, in
Figure~\ref{pentagon}.

A dual situation arises when there are exactly two edges incident to any
vertex in ${\mathcal V}(n)$, that is, the basis elements $v_i$ appear in exactly
two defining relations. In this case the same procedure of concatenating
edges, deleting vertices, and giving the new edge the label of the
deleted vertex creates a graph with edges labeled
by ${\mathcal V}(n)$ and vertices labeled by ${\mathcal V}'(n)$. In this
case we say that the defining relations are {\em dual graphlike}
and call the associated graph $\G'(\T(n))$ the
{\em dual Tel-A-graph}. The defining relations for the operad 
{\it Lie} and $n=4$  are dual graphlike, as shown in the next example.
The existence of an orientation for the dual Tel-A-graph means that
it is possible to find a basis for $\F(E)(n)$ and for 
$J(n)=\free{\F(E)}R(n)$ such that any basis element for $\F(E)(n)$ 
which appears in one of the relations $J(n)$ appears in two relations 
once with a coefficient $+1$ and once with a coefficient $-1$.
\begin{figure}[ht]
\unitlength 1.85mm
\def\aa{{\bf a}}
\def\bb{{\bf b}}
\def\cc{{\bf c}}
\def\dd{{\bf d}}
\def\ee{{\bf e}}
\def\ff{{\bf f}}
\def\gg{{\bf g}}
\def\hh{{\bf h}}
\def\ii{{\bf i}}
\def\jj{{\bf j}}
\def\kk{{\bf k}}
\def\ll{{\bf l}}
\def\mm{{\bf m}}
\def\nn{{\bf n}}
\def\oo{{\bf o}}
\def\I{{\bf 1}}
\def\II{{\bf 2}}
\def\III{{\bf 3}}
\def\IV{{\bf 4}}
\def\V{{\bf 5}}
\def\VI{{\bf 6}}
\def\VII{{\bf 7}}
\def\VIII{{\bf 8}}
\def\IX{{\bf 9}}
\def\X{{\bf 10}}

\begin{center}
\begin{picture}(62.83,40.67)(9,7)
\thicklines

\put(30.00,10.00){\line(1,0){20.00}}
\put(50.00,10.00){\line(1,2){10.00}}
\put(60.00,30.00){\line(-4,3){20.00}}
\put(40.00,45.00){\line(-4,-3){20.00}}
\put(20.00,30.00){\line(1,-2){10.00}}
\put(20.00,30.00){\line(1,0){40.00}}
\put(30.00,10.00){\line(2,3){5.00}}
\put(35.00,17.50){\line(1,4){5.13}}
\put(40.00,38.00){\line(1,-4){5.13}}
\put(45.25,17.50){\line(2,-3){5.00}}
\put(35.00,17.50){\line(6,5){15.00}}
\put(45.17,17.67){\line(-5,4){15.33}}
\put(40.00,45.00){\line(0,-1){7.00}}
\put(40.00,45.00){\makebox(0,0)[cc]{$\bullet$}}
\put(20.10,30.00){\makebox(0,0)[cc]{$\bullet$}}
\put(30.00,10.00){\makebox(0,0)[cc]{$\bullet$}}
\put(50.20,10.00){\makebox(0,0)[cc]{$\bullet$}}
\put(59.83,30.00){\makebox(0,0)[cc]{$\bullet$}}
\put(40.10,38.00){\makebox(0,0)[cc]{$\bullet$}}
\put(30.00,30.00){\makebox(0,0)[cc]{$\bullet$}}
\put(50.00,30.00){\makebox(0,0)[cc]{$\bullet$}}
\put(45.10,17.60){\makebox(0,0)[cc]{$\bullet$}}
\put(35.00,17.50){\makebox(0,0)[cc]{$\bullet$}}
\put(19.17,30.00){\makebox(0,0)[rc]{$\V$}}
\put(29.17,39.00){\makebox(0,0)[rc]{$\aa$}}
\put(24.17,19.00){\makebox(0,0)[rc]{$\ee$}}
\put(61.00,30.00){\makebox(0,0)[lc]{$\II$}}
\put(57.00,21.00){\makebox(0,0)[lc]{$\cc$}}
\put(55.00,31.00){\makebox(0,0)[lc]{$\gg$}}
\put(40.00,47.00){\makebox(0,0)[cb]{$\I$}}
\put(41.00,40.67){\makebox(0,0)[cb]{$\ll$}}
\put(50.50,39.00){\makebox(0,0)[cb]{$\bb$}}
\put(29.00,9.33){\makebox(0,0)[rt]{$\IV$}}
\put(51.00,9.33){\makebox(0,0)[lt]{$\III$}}
\put(40.00,9.00){\makebox(0,0)[lt]{$\dd$}}
\put(41.00,37.83){\makebox(0,0)[lc]{$\IX$}}
\put(42.50,32.83){\makebox(0,0)[lc]{$\jj$}}
\put(50.50,29.50){\makebox(0,0)[lt]{$\VIII$}}
\put(45.50,25.50){\makebox(0,0)[lt]{$\hh$}}
\put(29.50,29.50){\makebox(0,0)[rt]{$\VI$}}
\put(34.00,31.00){\makebox(0,0)[rb]{$\ii$}}
\put(32.00,27.50){\makebox(0,0)[rt]{$\nn$}}
\put(26.00,31.00){\makebox(0,0)[rb]{$\mm$}}
\put(46.00,18.00){\makebox(0,0)[lb]{$\VII$}}
\put(48.30,14.00){\makebox(0,0)[lb]{$\ff$}}
\put(36.00,17.33){\makebox(0,0)[lt]{$\X$}}
\put(33.50,14.33){\makebox(0,0)[lt]{$\oo$}}
\put(34.00,23.33){\makebox(0,0)[lt]{$\kk$}}
\end{picture}\caption{\label{peterson}
The Lie-hedron.}
\end{center}
\end{figure}
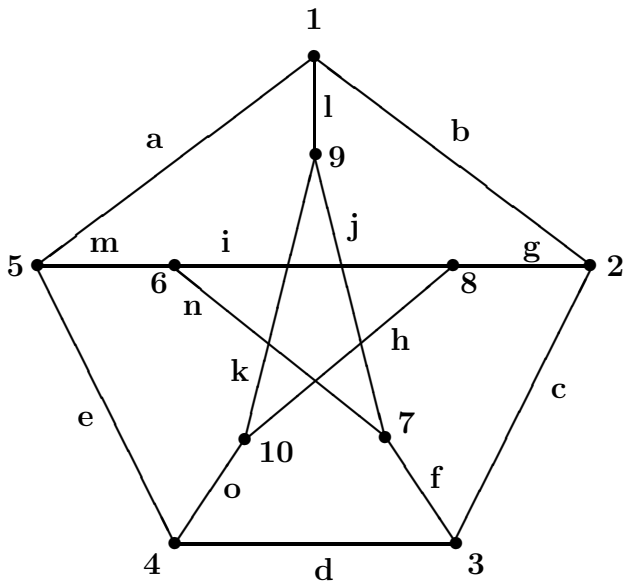

\begin{example}{\rm\
\label{zajicek_usacek}
For the Lie algebra operad $\Lie$ introduced in Example~\ref{Lie}
there is a dual Tel-A-graph constructed from the bipartite 
graph ${\mathcal T}(4)$ 
with vertices $({\mathcal V}(4),{\mathcal V}'(4))$ 
as given in Figure~\ref{matrix}. 
The resulting graph is  the famous Peterson graph, see
Figure~\ref{peterson}.

Let us describe $\T(3)$. Let 
$$ v_1 := [1[23]],\
v_2 := [2[31]],\
v_3 := [3[12]],\mbox{ and } r=\iota:=v_1+v_2+v_2. 
$$
Then  ${\mathcal V}(3) := \{v_1,v_2,v_3\}$ and
${\mathcal V}'(3) :=  \{r\}$. The bipartite graph $\T(3)$ 
has one trivalent vertex
$r$ and three univalent vertices $v_i$:
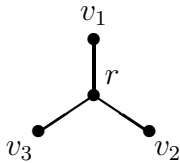
\begin{center}
{
\thicklines
\unitlength=.7pt
\begin{picture}(80.00,80.00)(0.00,0.00)
\put(0.00,0.00){\makebox(0.00,0.00){$v_3$}}
\put(10.00,10.0){\makebox(0.00,0.00){$\bullet$}}
\put(80.00,0.00){\makebox(0.00,0.00){$v_2$}}
\put(70.00,10.0){\makebox(0.00,0.00){$\bullet$}}
\put(40.00,73.00){\makebox(0.00,0.00){$v_1$}}
\put(40.00,60.00){\makebox(0.00,0.00){$\bullet$}}
\put(50.00,40.00){\makebox(0.00,0.00){$r$}}
\put(40.00,30.00){\makebox(0.00,0.00){$\bullet$}}
\put(40.00,30.00){\line(3,-2){30.00}}
\put(40.00,30.00){\line(-3,-2){30.00}}
\put(40.00,60.00){\line(0,-1){30.00}}
\end{picture}}
\end{center}
so it is neither graphlike, nor dual 
graphlike. In fact, $n = 4$ is the only value for which $\T(n)$ is (dual) 
graphlike.\endex}\end{example}

There is a particular class of 
operads whose defining relations are always graph-like.
Such operads $\P\in$ \kVect\ have  a 
presentation $\P=\prez ER$ in which we have chosen 
a basis $B(n)$ for $E(n)$ and a basis $K(n)$ for the relations 
$R(n)$ such that $\pi(r_i)=p_i-q_i$ for each $r_i\in K(n)$. 
For this class we can establish 
the relation between operad theory and
categorical coherence mentioned in the introduction. 
The following formulation arose from discussions 
between the first author and Tom~Fox.

If $\L$ is an operad in the monoidal category $\Sets$ of sets,
we can form the collection 
$$
\P={\rm Span}_{{\bold k}}(\L) =
\{{\rm Span}_{{\bold k}}(\L)(n)\}_{n\geq 1},
$$
where ${\rm Span}_{{\bold k}}(\L)(n):= {\rm Span}_{{\bold k}}(\L(n))$
is the $\bk$-vector space spanned by $\L(n)$. This collection has an 
obvious {\bf k}-{\bf Vect}-operad structure induced by the 
\Sets-operad structure of $\L$. If in addition each $\L(n)$ is
graded, then the $\bk$-vector space spanned by $\L(n)$ has a \kVect-operad 
structure. The preferred basis for each $\P(n)$ is given by the elements of
$\L(n)$. In short, an operad has a presentation with defining relations
of the type $\pi(r_i)=p_i-q_i$ if and only if it is a  
$\bk$-linearization of an  operad defined over the category of sets. 

The operad $\L$ can presented 
as $\L=\F_{\bf S}/\sim$, for $\F_{\bf S}$ a 
free \Sets-operad and $\sim$ an equivalence relation
given by a list of couples  
$p_i\sim q_i$ for $p_i,q_i\in \F_{\bf S}(B)(n)$. 
Each identification corresponds to a commutative 
diagram for the maps describing the corresponding $\P$-algebra. 

To make the exposition easier, in speaking about operads defined by
commutative diagrams, we shall restrict to non-$\Sigma$ operads.
The formal definition is  

\begin{definition}
\label{osta}
Let $\P$ be a \kVect\ operad. We say that $\P$ is 
{\em defined by commutative diagrams} 
if there exists a \Sets-operad $\L$ and an isomorphism
\begin{equation}
\label{/}
{\rm Span}_{{\bold k}}(\L) \cong \P
\end{equation}
of \kVect-operads.
\end{definition}

Categorical coherence means the commutativity of a family of  diagrams. 
Being defined by commutative diagrams is a crucial property of an
operad which
allows us to establish a relation between the theory of operads and the
theory of categorical structures. The following proposition
follows immediately from the previous discussion.

\begin{proposition}
\label{;}
If $\P$ is an operad defined by commutative diagrams, 
then for any $n\geq 3$, the bipartite graphs $\T_\P(n)$ 
are graphlike, and the corresponding Tel-A-graphs are oriented.
\end{proposition}
\noindent
The next proposition relates the Tel-A-graph $\G$ and the space 
${\mathcal D}$ of coherence relations.
 
\begin{proposition}
\label{kuba}
Let $\P\cong{\rm Span}_{{\bold k}}(\L)$ be an operad defined by 
commutative diagrams, where $\L$ is the corresponding
\Sets-operad presented by a collection
of generators $B$ and a collection of relations $K\in \F_{\bf S}(B)$. Let
$\G(n) = \G(\T(n))$ for $n\geq 3$ be the corresponding oriented Tel-A-graphs. 
We interpret the graph $\G(n)$ as a 1-dimensional simplicial complex. 
To each closed oriented path $\gamma\in\G(n)$ we can associate the
a coherence relation. This correspondence defines a natural isomorphism 
of ${\bf k}$-vector spaces:
\begin{equation}
\label{lemon}
H_1(\G(n), {\bf k})\cong \D(n).
\end{equation}
\end{proposition}

\begin{proof}
Let ${\mathcal V}(n)=\setdots{v_1}{v_b}=\F_{\bf S}(B)(n)$ and
${\mathcal V}'(n)=\setdots{r_1}{r_s}$.
For each $l\leq j \leq s$, we have $\pi(r_j) =
v_{a_j}- v_{b_j}$ for some $1\leq a_j,b_j\leq b$. We can orient the graph
$\G(n)$ so that $\partial (r_j)=v_{a_j}- v_{b_j},$ where $\partial$ is the
boundary operator for the one dimensional oriented simplicial complex
$\G(n)$ and then $\sum a_i r_i\in\Ker\partial=\Ker\pi $ if and only if 
$\sum a_i r_i$ expresses a coherence relation.
Moreover $\Ker\partial  \cong H_1(\G(n), {\bf k})$ 
since there are no boundaries in dimension $1$.
\end{proof}

The diagrams whose commutativity is required for
coherence correspond to elements of the
of the edge-path fundamental group of $\G(n)$, and therefore are
more closely related to $H_1(\G(n), {\mathbb Z})$ 
than to $H_1(\G(n), {\bf k})$. 
Therefore it is natural to consider integral homology. 
The fact that there are no $1$-cycles which are boundaries in $\G(n)$
implies that there is no torsion in $H_1(\G(n), {\mathbb Z})$, so
$$H_1(\G(n), {\bf k})\cong H_1(\G(n), {\mathbb Z})\otimes_{\mathbb Z}\k.$$ 

We also define an integral form for the operad $\L$.

\begin{definition}
Given the \Sets\ operad $\L=\F_{\bf S}(B)/\sim$ with $\sim$ consisting of a 
list of couples $p_{n,i}\sim p'_{n,i}$ for 
$p_{n,i},p'_{n,i}\in \F_{\bf S}(B)(n)$, we
define an operad $\L_{\mathbb Z}:=\prez {E_{\mathbb Z}}{R_{\mathbb Z}}$ in the 
symmetric monoidal category \AbGrp\ 
of abelian groups. This operad is generated by the collection of free
abelian groups 
$E_{\mathbb Z}(n):={\mathbb Z}^{B(n)}$ and the collection of relations 
$R_{\mathbb Z}(n):=\bigoplus_i {\mathbb Z}\cdot(p_{n,i}-p'_{n,i})$.
\end{definition}

The collection of coherence relations $\D_{\L_{\mathbb Z}}$ of the operad 
$\L_{\mathbb Z}$ is the $\F(E_{\mathbb Z})$-module defined completely
analogously to $\D_\P$ in Definition~\ref{61} and the coherence constraints 
$\C_{\L_{\mathbb Z}}$ is the corresponding collection of indecomposables. 
The following proposition describes the obvious relation between 
coherence data for $\P ={\rm Span}_{{\bold k}}(\L)$ 
and $\L_{\mathbb Z}$.

\begin{proposition}\label{Zbasis}
Let $\P$ be an operad defined by commutative diagrams, with associated
sets operad $\L=\prez BK$ such that there is an isomorphism (\ref{/}). Then
$$\D_\P\cong \D_{\L_{\mathbb Z}}\otimes_{\mathbb Z}\bfk,\quad \C_\P\cong 
\C_{\L_{\mathbb Z}}\otimes_{\mathbb Z}\bfk,$$
and 
\begin{equation}
\label{den_pred_odletem}
H_1({{\sf G}}_\P(n),{\mathbb Z})\cong \D_{\L_{\mathbb Z}}.
\end{equation}
\end{proposition}

Before ending this section we make a final remark
about operads with defining relations which are dual graphlike
as in the Lie-hedron. Recall that in the dual Tel-A-graph the vertices 
correspond to relations and the edges correspond to basis elements 
appearing in the relation.  In this case  the coherence relations 
do not correspond to the space of $1$-cycles, but to certain $0$-chains.
If it is possible to orient the dual Tel-A-graph $\G'(\T(n))$,
then the sum of all the vertices in a connected component defines a 
coherence constraint since it represents a
linear combination of relations such that each basis element (edge)
appears once with a plus sign and once with a negative sign. 
It is easy to see from the matrix in Figure~\ref{matrix} that by changing the 
sign of four of the relations, defining
${\bf i'}=-{\bf i}$ for ${\bf i}\in \{{\bf 3,6,7,9}\}$,
 we get a matrix with each column having two nonzero entries,
one $+1$ one $-1$, which allows us to orient the Lie-hedron.
Then the coherence relation
$$
-\ell:={\bf 1}+{\bf 2} +{\bf 3'} 
+{\bf 4}+{\bf 5} +{\bf 6'} +{\bf 7'} +{\bf 8} +{\bf 9'} +{\bf 10}
$$
is a just a sum over the vertices. The fact that the Lie-hedron is
connected implies that the space of coherence constraints in $n=4$
is one dimensional. We shall discuss this issue further in 
Example~\ref{modra_propiska}. 

The following two useful formulas 
(one for the non-$\Sigma$, one for the symmetric
case) compute the size of the sets
${\mathcal V}(4)$ and ${\mathcal V}'(4)$ for a 
quadratic operad $\P =\prez ER$. The formulas are immediate 
consequences of the description of the free operad and free operad 
module in terms of trees.
\begin{eqnarray}
\label{pkc}
\dim({\mathcal V}(4))\!=\! 5\dim(E)^3\!,\ 
\dim({\mathcal V}'(4))\!=\! 5 \dim(E)\dim(R) \hskip -.5cm &&\mbox{ (the
non-$\Sigma$ case)}
\\
\label{pks}
\dim({\mathcal V}(4))\!=\! 15\dim(E)^3\!,\, 
\dim({\mathcal V}'(4))\!=\! 10 \dim(E)\dim(R)\hskip -.7cm &&\mbox{
(the
symmetric case)}
\end{eqnarray}
Summing up the above results, we see for a Koszul operad $\P$
the following nice formula which relates the topological properties
of the graph $\G(4)$ and the algebraic properties of the operad $\P$:
\[
\dim(\C_\P)= \dim(\P^!(4)) = \dim(H_1(\G(4))).
\]

\section{Relation to Mac Lane coherence}
\label{5}

In this section we show how our theory gives coherence theorems \`a~la Mac~Lane. 
We will deal with operads given by commutative diagrams, in the sense of 
Definition~\ref{osta}.
For these operads, the defining relations are always graphlike
(Proposition~\ref{;}) with $\G(n)$ the corresponding
Tel-A-graph in degree $n$.

Let $\P$ be an operad defined by commutative diagrams, and assume that
the isomorphism (\ref{/}) has been fixed once and for all, as well as the 
presentation $\L=\prez BK$. Assume we have oriented the 
Tel-A-graphs $\G(n)$ as in the proof of Proposition~\ref{kuba} so
that $\pi(r)=\partial(r)$. Define the functions $\eta,\xi : {\mathcal V}'(n)
\to {\mathcal V}(n)$ by
$\pi(r)=\eta(r)-\xi(r)$ for any edge $r$ of $\G(n)$.

\begin{definition}
\label{America_West}
A $(B,K)$-structure on a category ${\mathcal C}$ consists  of:
\begin{itemize}
\item[(i)]
 An assignment of an $n$-polyfunctor 
$\Phi(b)$ for each element $b\in B(n)$ and thus, 
by functorial extension, an 
$n$-multifunctor  $\Phi(p)$ for every element $p\in \F_{\bf S}(B)(n)$.
\item[(ii)]
An assignment of a natural isomorphism $a(r)$ 
between $\Phi(\xi (r))$ and $\Phi(\eta(r))$ for each $r\in K$ and
thus, by functorial extension, for every element
$r\in \F_{\bf S}(B)\langle K\rangle$.
\end{itemize}
\end{definition}

The $(B,K)$-structure is {\em coherent} if and only if whenever there is 
a composition of the natural transformations 
$a(r)$ and their inverses connecting a given pair of 
multi-functors $\Phi(p)$ and $\Phi(q)$, 
then all possible compositions connecting
these two multi-functors give the same natural transformation. 
This is equivalent 
to saying that whenever there is a composition
of the $a(r)$'s and their inverses which connects the 
multi-functor $\Phi(p)$ to itself, the composition
is the identity transformation.

Instead of a $(B,K)$-structure we will often speak simply about a
$\calP$-structure, $\calP = \Span_{\bf k}(\calL)$. The  $(B,K)$-%
notation, on the other hand, underlines the fact that the structure of
Definition~\ref{America_West} is very explicitly related to the
prezentation $\prez BK$ of $\calL$.
 
Each closed path 
$\gamma \in \G(n)$ determines a diagram $D(\gamma)$ with arrows 
corresponding to natural isomorphisms of functors. 
By definition, the $(B,K)$-structure on the category 
${\mathcal C}$ is coherent if and only if all these diagrams are commutative.

\begin{proposition}
\label{23}
Let $n\geq 3$ and let $\setdots{\gamma_1}{\gamma_d}$ be a sequence of
closed paths of the graph $\G(n)$ forming a basis of
$H_1(\G(n),{\mathbb Z})$. Then the diagram of natural isomorphisms
$D(\gamma)$ is commutative for any
closed path $\gamma$ of $\G(n)$ if and only if it is commutative for
any $\gamma \in \setdots{\gamma_1}{\gamma_d}$.
\end{proposition}

\begin{proof} 
Since the graph is an one-dimensional simplicial complex, the fundamental 
group is isomorphic to the free group on a set of closed paths. The  set of 
generators also forms a homology basis. 
Conversely, each set of closed paths that forms a basis for the
cohomology generates the fundamental group.

Thus in order to guarantee 
commutativity of all diagrams $D(\gamma)$ 
it is enough to check that the diagrams $D(\gamma_i)$, $1\leq i \leq d$, 
are commutative.
\end{proof}

As a corollary we formulate the following `classical'
coherence theorem which involves the coherence constraints 
$\C_{\L_{\mathbb Z}}$  of the operad $\L_{\mathbb Z}$ in the
category \AbGrp. 

\begin{theorem}
\label{21}
A necessary and sufficient condition for the coherence of a 
$(B,K)$-structure on the category ${\mathcal C}$ is that the diagrams
$D(\gamma)$ are commutative for a set ${\mathcal B}$ 
of closed paths $\gamma$ such that
the images \{$p(\pi(\gamma));\ \gamma \in {\mathcal B}\}
$ form a basis for ${\C}_{\L_{\mathbb Z}}$. 

Thus, in the case of a Koszul operad there are
$\dim(\P^!(4))$-diagrams whose commutativity is necessary and sufficient for
coherence.
\end{theorem}

\begin{proof} 
The commutativity of all diagrams $D(\gamma)$, $\gamma \in {\mathcal B}$, 
is clearly necessary for the coherence. We need to prove that it
is sufficient.

By Proposition~\ref{23}, the condition for 
coherence is the commutativity of the diagrams $D(\gamma)$ for a basis of 
$H_1(\G(n),{\mathbb Z})$. 
By assumption the diagrams with $p(\pi(\gamma))$ 
forming a basis for ${\C}_{\L_{\mathbb Z}}$ 
are commutative. By~(\ref{den_pred_odletem}) of
Proposition~\ref{Zbasis},
these together with the `obvious relations' give
a ${\mathbb Z}$-basis for $H_1(\G(n),{\mathbb Z})$. 

Thus, to finish the proof of the first part of the theorem, it is
enough to show that each `obvious relation' corresponds to a
commutative diagram.
The first type of the `obvious relation'  
\[
x(b_1,\ldots,b_{s-1},\pi(y),b_{s+1},\ldots,b_l) -
\pi(x)(b_1,\ldots,b_{s-1},y,b_{s+1},\ldots,b_l) =0
\]
is, for $b_1,\ldots,b_{s-1},b_{s+1},\ldots,b_l \in \F_{\bf S}(B)$, 
$x,y \in \F_{\bf S}(B)\langle K
\rangle$,
represented by the simple closed path $\omega$ (recall $\pi(x) =
\eta(x)-\xi(x)$, $\pi(y) = \eta(y) - \xi(y)$):

\vglue3mm
\def\maletecky{{\hskip -.4mm .\hskip -.2mm .\hskip -.2mm.\hskip -.2mm}}
\def\pomoc#1#2{{#1}(x)(b_1,\maletecky,b_{s-1},%
               {#2}(y),b_{s+1},\maletecky,b_l)}
\def\Pomoc#1#2{{#1}(b_1,\maletecky,b_{s-1},%
               {#2},b_{s+1},\maletecky,b_l)}

\vglue .55cm
\begin{center}
\square{\pomoc\xi\xi \hskip 4.3cm}{\hskip 4.3cm\pomoc\eta\xi}%
       {\pomoc\xi\eta \hskip 4.3cm}{\hskip 4.3cm\pomoc\eta\eta}%
       {\Pomoc x{\xi(y)}}{\Pomoc {\xi(x)}y \hskip -4mm}%
       { \hskip -4mm\Pomoc{\eta(x)}y}{\Pomoc x{\eta(y)}}
\end{center}

\vglue .8cm
\noindent 
for which $D(\omega)$ is
commutative by naturality. Also the
second `obvious' relation
\begin{eqnarray*}
\lefteqn{
b(b_1,\ldots,b_{s-1},\pi(x),b_{s+1},
\ldots,b_{t-1},y,b_{t+1},\ldots,b_l) -\hskip 2cm}
\\
&&\hskip 2cm -b(b_1,\ldots,b_{s-1},x,b_{s+1},
\ldots,b_{t-1},\pi(y),b_{t+1},\ldots,b_l)=0
\end{eqnarray*}
can be represented by a simple closed 
path $\omega$ such that $D(\omega)$ is commutative in 
an obvious similar way, again by the naturality.

The second part of the theorem follows from the description of the
coherence constraints
$\C_{\mathcal P}$ given in Theorem~\ref{main}.
\end{proof}

\begin{example}{\rm
\label{25}
(continuation of Example~\ref{hedges})
The non-$\Sigma$ operad $\Ass$ for associative algebras is 
quadratic Koszul and
$\Ass = \Ass^!$, therefore $\C_\Ass = \Ass(4) = \Span(p)$, where $p$
{\em must\/} be the element of~(\ref{assrelate}) --
there is no other choice! We get from~(\ref{pkc}) that $\dim({\mathcal V}(4))=
\dim({\mathcal V}'(4))= 5$, and $\G(4)$ is the pentagon.
An $\Ass$-structure on a category ${\mathcal C}$ is the same as a
multiplication on ${\mathcal C}$ with an associativity isomorphism as it
was discussed in the introduction. The
coherence of $({\mathcal C},\BOX,a)$ in the sense of the above
definitions coincides with Mac~Lane's definition,
and Theorem~\ref{21} gives Mac~Lane's
celebrated coherence result.
\endex}\end{example}

\begin{example}{\rm\
\label{[]}
Consider the algebraic structure consisting of a
vector space $V$ and two bilinear maps $\circ,\bullet : V\ot V\to V$ which
satisfy
\begin{eqnarray*}
a\circ(b\circ c) = (a\circ b) \circ c, && 
a\circ (b \bullet c) = (a\circ b) \bullet c,
\\
a\bullet (b\circ c) = (a\bullet b) \circ c, 
&& a\bullet (b\bullet c) = (a\bullet b) \bullet c.
\end{eqnarray*}
These algebras were introduced in~\cite{markl:zebrulka} and called
{\em nonsymmetric Poisson algebras\/}. The corresponding operad
$\K$ is Koszul, quadratic self-dual
(see again~\cite{markl:zebrulka}) and it 
has a quadratic presentation $\K = \prez ER$ with $\dim(E)= 2$ and
$\dim(R)= 4$. We easily calculate that $\dim(\K(4))= 8$, thus, by
Theorem~\ref{main},  
$\dim(\C_\K) = \dim (\K^!(4)) =\dim(\K(4))= 8$, with a basis
consisting of eight pentagons shown in Figure~\ref{eight}.
\begin{figure}[ht]
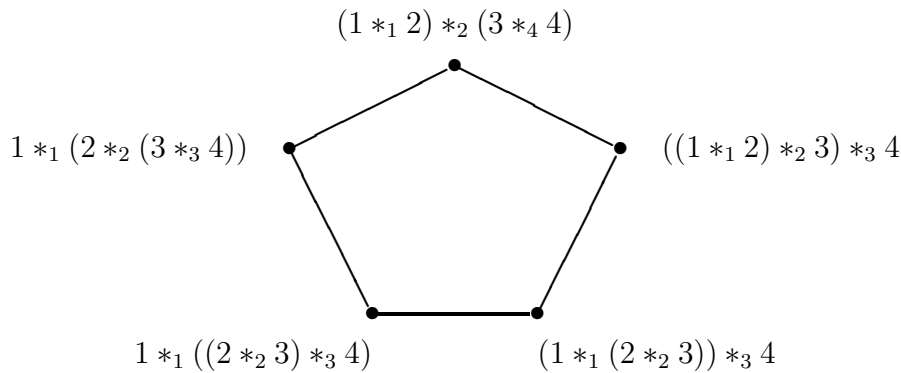

\setlength{\unitlength}{1.1cm}
\pentagonline{(1*_1(2*_23))*_34}{((1*_12)*_23)*_34}%
{(1*_12)*_2(3*_44)}{1*_1(2*_2(3*_34))}%
{1*_1((2*_{2} 3)*_3 4)}
\caption{
\label{eight}
The graph $\G(4)$ for nonsymmetric Poisson algebras.
The triple $(*_1,*_2,*_3)$
runs through all eight possible combinations $(\circ,\circ,\circ)$,
$(\bullet,\circ,\circ)$, $(\circ,\bullet,\circ)$, $(\circ,\circ,\bullet)$, $(\circ,\bullet,\bullet)$,
$(\bullet,\circ,\bullet)$, $(\bullet,\bullet,\circ)$ and $(\bullet,\bullet,\bullet)$.}
\end{figure}

A $\K$-structure on
${\mathcal C}$ consists of two covariant bifunctors
$\BOX_1,\BOX_2 :
{\mathcal C}\times {\mathcal C} \to {\mathcal C}$ and four natural
transformations
\begin{eqnarray*}
a_{11} : \BOX_1(\id \times \BOX_1) \to \BOX_1(\BOX_1 \times \id),
&&
a_{12}:\BOX_1(\id \times \BOX_2) \to \BOX_2(\BOX_1 \times \id),
\\
a_{21}:\BOX_2(\id \times \BOX_1) \to \BOX_1(\BOX_2 \times \id),
&&
a_{22}:\BOX_2(\id \times \BOX_2) \to \BOX_2(\BOX_2 \times \id).
\end{eqnarray*}
This structure is a variant of a {\em weakly distributive
category\/}~\cite[Definition~1.1]%
{blute-cockett-seely-trimble:JPAA96}.
These categories are important for linear logic; we intend to
discuss the applications of our theory to this direction in another
paper.
By Theorem~\ref{21}, such a $\K$-structure is coherent if and only if
the eight pentagonal diagrams corresponding to those in Figure~\ref{eight}
commute.\endex}\end{example}

\begin{example}{
\rm\
\label{digebres}
Let us discuss the following bizarre objects introduced by
Loday
in~\cite{loday:digebres}. The importance of this example is that
some coherence constraints
{\em will not be the pentagons\/}.

By a {\em digebra\/} we mean a vector space
$V$ together with two bilinear operations $\circ$ and $\bullet$ which
satisfy the following axioms:
\begin{eqnarray*}
&x \circ (y\circ z) = (x\circ y)\circ z = x \circ (y\bullet z),&
\\
&(x\bullet y)\circ z = x\bullet (y\circ z),&
\\
&(x\circ y)\bullet z = x\bullet (y \bullet z) = (x\bullet y)\bullet z.&
\end{eqnarray*}
Let $\D$ be the corresponding non-$\Sigma$ operad. It has a
quadratic presentation $\D = \prez ER$ with $\dim(E)= 2$ and
$\dim(R)= 5$, and formula~(\ref{pks}) gives that $\dim(\V(4))= 40$
and $\dim(\V'(4))= 50$. 

As it was proven in~\cite{loday:digebres},
the
operad $\D$ is Koszul. It is  easy to compute that $\dim(\D(n))= n$,
and Proposition~\ref{za} says that $\dim(\C_\D)=14$. The graph
$\G(4)$ is complicated (it has 40 vertices and 50 edges),
but we know, by Proposition~\ref{23},
that there exist 14 closed cycles in $\G(4)$
which generate the coherence constraints. These cycles are shown on
Figure~\ref{24}.

\begin{figure}
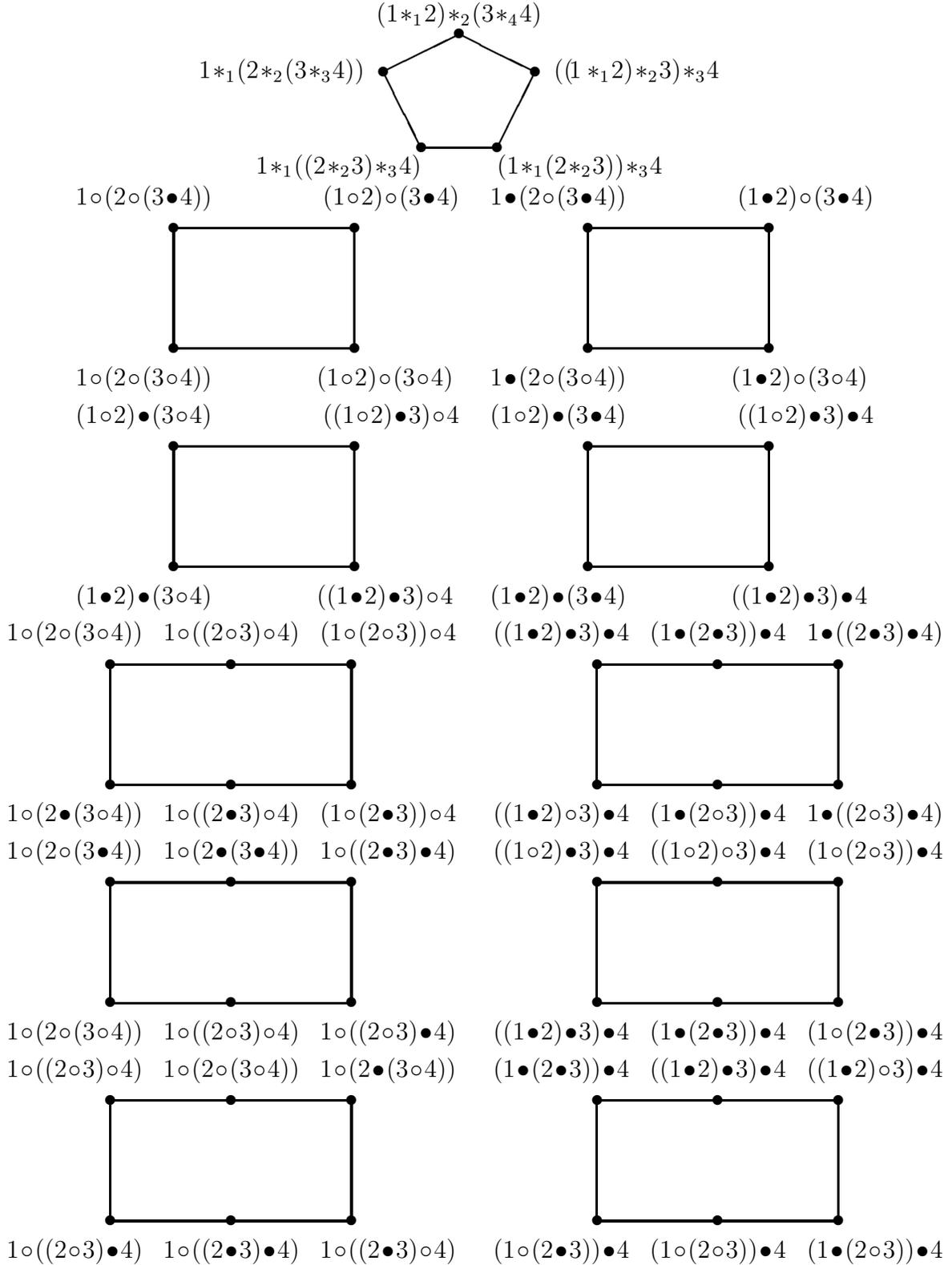

\setlength{\unitlength}{.63cm}
\pentagonline{(1\!*_1\!(2\!*_2\!3))\!*_3\!4}%
{((\!1*_1\!2)\!*_2\!3)\!*_3\!4}%
{(1\!*_1\!2)\!*_2\!(3\!*_4\!4)}{1\!*_1\!(2\!*_2\!(3\!*_3\!4))}%
{1\!*_1\!((2\!*_{2}\! 3)\!*_3\! 4)}

\vskip4mm
\begin{center}
\ssquare{\ix{\CC}{\CC}{\BB}}{\ii{\CC}{\CC}{\BB}}%
 {\ix{\CC}{\CC}{\CC}}{\ii{\CC}{\CC}{\CC}}
\hskip2cm
\ssquare{\ix{\BB}{\CC}{\BB}}{\ii{\BB}{\CC}{\BB}}%
 {\ix{\BB}{\CC}{\CC}}{\ii{\BB}{\CC}{\CC}}
\end{center}

\begin{center}
\ssquare{\ii{\CC}{\BB}{\CC}}{\iii{\CC}{\BB}{\CC}}%
 {\ii{\BB}{\BB}{\CC}}{\iii{\BB}{\BB}{\CC}}
\hskip2cm
\ssquare{\ii{\CC}{\BB}{\BB}}{\iii{\CC}{\BB}{\BB}}%
 {\ii{\BB}{\BB}{\BB}}{\iii{\BB}{\BB}{\BB}}
\end{center}

\begin{center}
\hex{\ix{\CC}{\CC}{\CC}}{\vv{\CC}{\CC}{\CC}}{\iv{\CC}{\CC}{\CC}}%
 {\ix{\CC}{\BB}{\CC}}{\vv{\CC}{\BB}{\CC}}{\iv{\CC}{\BB}{\CC}}
\hskip2.3cm
\hex{\iii{\BB}{\BB}{\BB}}{\iv{\BB}{\BB}{\BB}}{\vv{\BB}{\BB}{\BB}}%
 {\iii{\BB}{\CC}{\BB}}{\iv{\BB}{\CC}{\BB}}{\vv{\BB}{\CC}{\BB}}
\end{center}

\begin{center}
\hex{\ix{\CC}{\CC}{\BB}}{\ix{\CC}{\BB}{\BB}}{\vv{\CC}{\BB}{\BB}}%
 {\ix{\CC}{\CC}{\CC}}{\vv{\CC}{\CC}{\CC}}{\vv{\CC}{\CC}{\BB}}
\hskip2.3cm
\hex{\iii{\CC}{\BB}{\BB}}{\iii{\CC}{\CC}{\BB}}{\iv{\CC}{\CC}{\BB}}%
 {\iii{\BB}{\BB}{\BB}}{\iv{\BB}{\BB}{\BB}}{\iv{\CC}{\BB}{\BB}}
\end{center}

\begin{center}
\hex{\vv{\CC}{\CC}{\CC}}{\ix{\CC}{\CC}{\CC}}{\ix{\CC}{\BB}{\CC}}%
 {\vv{\CC}{\CC}{\BB}}{\vv{\CC}{\BB}{\BB}}{\vv{\CC}{\BB}{\CC}}
\hskip2.3cm
\hex{\iv{\BB}{\BB}{\BB}}{\iii{\BB}{\BB}{\BB}}{\iii{\BB}{\CC}{\BB}}%
 {\iv{\CC}{\BB}{\BB}}{\iv{\CC}{\CC}{\BB}}{\iv{\BB}{\CC}{\BB}}
\end{center}

\vskip-4mm
\caption{\label{24}
The graph $\G(4)$ for digebras. It consists of four pentagons
(the triple $(*_1,*_2,*_3)$ runs through $(\circ,\circ,\circ)$,
$(\bullet,\circ,\circ)$, $(\bullet,\bullet,\circ)$ and $(\bullet,\bullet,\bullet)$), four squares
and six hexagons.}
\end{figure}
}\end{example}

\section{Operads and their algebras in a category of modules}
\label{6}

Recall~\cite{may:1972} 
that a $\P$-algebra structure on a (graded) $\bfk$-vector space $U$ is an
operad map $A: \P \to \End U$ from the operad $\P$ to the
endomorphism operad $\End U$ of $U$, where all structures are
considered in the category ${\tt Vect}_{\bf k}$ or \kVect\
with multiplication given by the (graded) tensor product over $\bk$, 
$\otimes:=\otimes_\bk$. If $U$ is a left $V$-module
over a  $\bk$-algebra $V$, then $\End U(n)$ has a natural
left $V$- \hskip .5mm 
right $V^{\ot n}$-module structure. This allows us to consider
generalized $\P$-algebras satisfying axioms that may involve coefficients
from $V$ and its tensor powers. In this section we study the way in which
the structure of the operad $\P$ imposes conditions on $V$ which are
necessary for there to be a reasonable concept of generalized  
$\P$-algebras. 

An example of this type was given in the introduction, where we considered
$\Phi$-associative algebras and derived the 
pentagon identity on the associator for the bialgebra $V$.  

Fix a unital, associative, local $\bfk$-algebra 
$V$, with augmentation $\eps:V\rightarrow \bfk$
(i.e., each $v \in V$ with $\eps(v)\not = 0$ is invertible) and an
operad $\P = \prez ER$.
In order to take tensor products of $V$-modules we must assume that, for each $n$
with $E(n)\not = 0$, we are given a $\bfk$-linear algebra homomorphism
(a `diagonal') $\Delta^n
: V \to V^{\ot n}$ such that  
\begin{equation}
\label{mysicka}
(\epsilon^{\ot n}) \Delta^n(v) = \epsilon(v),
\mbox { for each $v\in V$,}
\end{equation}
where we identify, in the left hand side of~(\ref{mysicka}),  
$\bfk$ with $\bigotimes_{\bfk}^{n}\bfk$.
When $\P$ is quadratic, then $\Delta = \Delta^2 : V \to V\ot V$ and
$(V,\cdot,\Delta,1,\eps)$ is an ordinary augmented unital associative, non
necessarily coassociative, bialgebra.

For $\P = \prez ER$ and  $n\geq 1$ we have
the basic exact sequence of $\boldk$-vector spaces, which
 defines $\P$ as a quotient of the free operad generated by $E$,
\begin{equation}
\label{52}
0
\longrightarrow \ker(\pi)(n)
\longrightarrow
\free{\F(E)}R(n)
\stackrel\pi\longrightarrow
\F(E)(n)
\longrightarrow \P(n)  \longrightarrow 0.
\end{equation}
This identifies
$\P(n)$ to $\Coker {\pi}{(n)}$.
We now define  an analog of this sequence which brings
$V$ into play and allows us to consider
generalized coherence conditions. Let us discuss the non-$\Sigma$
case first.

We already observed in Section~\ref{4} that elements of
$\F(E)(n)$ are represented by a sum of planar rooted
trees with vertices labeled by
elements of $E$. Such a labeled
tree $t$ determines uniquely a
bracketing $b = b_t$ of $n$ indeterminates or, equivalently,
an iterated comultiplication denoted $\Delta^{n,b}$.
Form the {\em $V$-relative free operad $\F_V(E)$\/} as follows:
\begin{equation}
\label{v-op}
\F_V(E)=\bigoplus_{n\geq 1} \F_V(E)(n),
\mbox{ with }
\F_V(E)(n):=\F(E)(n)\otimes_{\boldk} V^{\otimes n}.
\end{equation}
Each $\F_V(E)(n)$ has an obvious right $V^{\ot n}$-module structure.
We give it  a left $V$-module structure by defining
\begin{equation}
\label{leftV}
v\cdot(t\otimes \vec u)=t\otimes \Delta^{n,b_t}(v)\vec u,
\ \vec u = \orada{u_1}{u_n} \in V^{\ot n}.
\end{equation}
Given  $ t\otimes \vec u\in \F_V(E)(n)$, where
$\vec u  =u_{1}\otimes\ldots\otimes u_{n}$,
the `operadic' composition of this element  with the tensor
product of
$n$ elements $ t_i\otimes \vec v_i\in \F_V(E)(a_i)$,
$i=1,\ldots n$, is defined as
\begin{eqnarray}
\label{IIS}
\lefteqn{
\gamma((t\!\ot\! \vec u);(t_1\!\ot\! \vec v_1)\!\ot\! \cdots
\!\ot\!(t_n\!\ot\! \vec v_n)):=}
\\
\nonumber
&&:= t(t_1\!\ot\!\cdots\!\ot\!
t_n)\!\ot\!\Delta^{a_1,b_{t_1}}(u_{1}) \vec v_1\!\ot\!\ldots
\!\ot\! \Delta^{a_n,b_{t_n}} (u_{n})\vec v_n,
\end{eqnarray}
where $t(t_1\otimes\cdots\ot t_n)$ is the composition in $\F(E)$.
Heuristically, we can say that the composition moves
the interior coefficients $v_{i}$ across the tree
$t_i$ using the comultiplication $\Delta^{a_i, b_{t_i}}$.
Clearly this defines on $\F_V(E)$ the structure of a non-$\Sigma$
operad. In degree $n$ it is a free right $V^{\otimes n}$-module
on the
$\bk$-linear space $\F(E)(n)$.

In the symmetric case we define the right action of the symmetric
group on $\F_V(E)$ by
\[
[t\otimes (u_1\ot \cdots \ot u_n)]\sigma = t\sigma \ot
(u_{\sigma(1)}\ot \cdots \ot u_{\sigma(n)}),
\]
for $t\ot \vec u \in \F_V(E)(n)$ and $\sigma \in \Sigma_n$. 

To guarantee the  consistency of this action with the 
left $V$-module structure (\ref{leftV})
we define the comultiplication for the tree $t\sigma$ 
by
\begin{equation}
\label{treecomult}
\Delta^{n,b_{t\sigma}}(v):=\sigma^{-1}(\Delta^{n,b_t}(v)).
\end{equation}

Assume now that $R_V = \Coll{R_V}$ is a subcollection of 
$\F_V(E)$  such
that each $R_V(n)$ is left $V$- right $V^{\ot n}$-submodule.
 In the symmetric case we moreover require that $R_V(n)$
is $\Sigma_n$-closed. 

Then we form the operadic ideal 
$(R_V) \subset \F_V(E)$ and define 
the {\em $V$-relative operad $\P_V$\/} by 
\begin{equation}
\label{relative}
\P_V := \F_V(E)/(R_V).
\end{equation}
Observe that each $\P_V(n)$ is a left $V$- \hskip .5mm
right $V^{\ot n}$-module.

To formulate the concept of an operad algebra from this
point of view we use the usual endomorphism operad with
$\End U(n)=\Hom_{\boldk}(U^{\otimes n},U)$ considered
as a left $V$- and right $V^{\otimes n}$-module in the standard way:
\begin{eqnarray*}
(v\cdot \alpha)(u_1\otimes \cdots\otimes u_n)&=&
v\cdot(\alpha(u_1\otimes \cdots\otimes u_n),
\\
(\alpha\cdot(v_1\otimes \cdots\otimes v_n))
(u_1\otimes \cdots\otimes u_n)
&=&\alpha(v_1\cdot u_1\otimes \cdots\otimes v_n\cdot u_n),
\end{eqnarray*}
for $\alpha \in \End U(n)$, $\orada{u_1}{u_n}\in U^{\ot n}$, $v\in V$
and $\orada{v_1}{v_n}\in V^{\ot n}$.
A $\P_V$-algebra structure on the $V$-module
$U$ is given by a `$V$-relative'
operad map, i.e.~a
family of left $V$- right $V^{\ot n}$-module maps (equivariant, in the
$\Sigma$-case)
\[
A(n):\P_V(n)\longrightarrow \End U(n),\ n\geq 1,
\]
that are compatible with the operadic compositions.

\begin{example}{\rm\
\label{A}
Suppose that $U$ has a $\bk$-linear multiplication $*: U\ot U \to
U$ and assume that $(U,*)$ is a $V$-module algebra, that is,
\[
v(a*b) = v_{(1)}(a)* v_{(2)}b,
\mbox { for $v \in V$, $a,b \in U$,} 
\] 
where, as before, $v_{(1)}\otimes v_{(2)}$ stands for 
$\sum_iv_{(1)i}\otimes v_{(2)i}=\Delta(v)$ and
$\Delta$ is not necessarily coassociative.
We want to consider associativity in the category of $V$-modules.
Suppose that $U$ is a left $V = (V,\cdot)$-module and
replace $a*(b*c) = (a*b)*c$ by
\begin{equation}
\label{u2}
\Phi(a*(b*c)) = ((a*b)*c)
\end{equation}
where $\Phi\in V^{\ot 3}$ is an invertible element and
$\Phi(a*(b*c)) :=\sum (\Phi_1 a * (\Phi_2 b
*\Phi_3 c))$, if $\Phi = \sum \Phi_1\ot \Phi_2 \ot \Phi_3$.

These algebras are algebras over the $V$-relative operad $\Ass_\Phi$
defined as follows. 
In the notation introduced in Example~\ref{hedges}, 
let $R_V = R_{\Ass,\Phi}$ be a right $V^{\ot 3}$-submodule of
$\F_V(E)(3)$ generated by
\[
\tilde r_{{\Ass,\Phi}} = (1(23))\cdot \Phi - ((12)3) \in \F_V(E)(3).
\]
Then $\Ass_{\Phi} := \F_V(E)/(R_{\Ass,\Phi})$ 
is the $V$-relative 
operad describing algebras that satisfy~(\ref{u2}).%
\endex}\end{example}

The augmentation map $\epsilon : V \to \bfk$ induces, for each $n$, 
an augmentation (denoted by the same symbol) $\epsilon : V^{\ot n} \to
\bfk$. The right tensoring with $\bfk$ over $\Vn$ thus makes sense and
it defines  a
right-exact functor from the category of right $V^{\ot n}$-modules to
the category of $\bfk$-vector spaces. 
We denote this functor by $\eps(-)$. We also have, for any right
$V^{\ot n}$-module $M$, 
a canonical $\bfk$-linear epimorphism $\epsilon_M :
M \to \epsilon(M)$ given by $\eps_M(m) := m \ot_{\Vn}1 \in \eps(M)$. 

\begin{lemma}
\label{Nunyk}
For each $n\geq 1$, there is a canonical 
isomorphism 
\[
\epsilon(\F_V(E)(n)) \cong \F(E)(n).
\]
The canonical map 
$\epsilon_{\F_V(E)} : = 
\{\epsilon_{\F_V(E)(n)}:\F_V(E)(n)\rightarrow \F(E)(n)\}_{n\geq 1}$ 
is a homomorphism of~$\bfk$-operads.
\end{lemma}

\begin{proof}
The first statement of the lemma is obvious. The second part
immediately follows from definition~(\ref{IIS}) of the
operadic structure on $\F_V(E)$ and
compatibility~(\ref{mysicka}) of $\epsilon$ and $\Delta$'s.
\end{proof}

Let us introduce the central notion of this section.

\begin{definition}
\label{central_definition}
Assume that each $R_V(n)$ is a left $V$-submodule of $\F_V(E)(n)$,
is $\Sigma_n$-closed and is free as a right $V^{\ot n}$-module. 
Let $\P_V = \F_V(E)/(R_V)$ be as in~(\ref{relative}) 
and define $R := \epsilon(R_V)
\subset\F(E)$. Then we call the $V$-relative operad $\P_V$ a
$V$-relativization or $V$-quantization of 
the $\bfk$-operad $\P := \F(E)/(R)$. Let  
$\Set{\Rada {\tilde r^n}1{s(n)}} \subset \F_V(E)(n)$ be a
$V^{\ot n}$-basis of $R_V(n)$, $n\geq 3$. 
We require further that $({\Rada {r^n}1{s(n)}})$,
with $r^n_i := \epsilon(\tilde r^n_i)$, $1\leq i \leq s(n)$,  
are independent over
$\bk$ and thus form a basis for $R(n)=\epsilon(R_V(n))$. 
\end{definition}

Define the map $\pi_V: \F_V(E)\langle R \rangle \to \F_V(E)$ of
$\F_V(E)$-modules by $\pi_V(r_i) := \tilde r_i$, $1\leq i \leq s$.
We have the following $V$-relative analog
of~(\ref{52}):
\begin{equation}
\label{52q}
0  \longrightarrow  \ker(\pi_V)(n)
\longrightarrow
\F_V(E)\langle R\rangle(n)
\stackrel{\pi_V}{\longrightarrow}
\F_V(E)(n)
\longrightarrow \P_V(n)  \longrightarrow 0.
\end{equation}

Now apply the functor $\epsilon(-)$ on this sequence. Since clearly
$\epsilon(\F_V(E)\langle R \rangle (n)) \cong \F(E)\langle R \rangle
(n)$, $\epsilon(\F_V(E)(n)) \cong \F(E)(n)$ and $\epsilon(-)$ is right
exact, we obtain the exact sequence
\[
\free{\F(E)}R(n)
\stackrel\pi\longrightarrow
\F(E)(n)
\longrightarrow \epsilon(\P_V)  \longrightarrow 0
\]
where $\pi(r_i) = r_i$, $1\leq i \leq s$.
This implies the existence 
of the canonical
$\bfk$-isomorphism $\epsilon(\P_V(n))\cong\P(n)$. The
universal property of the kernel gives 
the canonical map
\begin{equation}
\label{canmap}
\rho:\epsilon(\Ker{\pi_V}(n))\to \Ker {\pi(n)}.
\end{equation}

Let $F:X\to Y$ be a map of two $\otexp{V}n$-modules and
let $f:\epsilon(X)\to \epsilon(Y)$ be a map of $\boldk$-vector
spaces. We say that $F:X\to Y$ is a {\em map over $f$\/} if
$\epsilon F=f\epsilon$.

\begin{definition}
\label{32}
Let $\P_V$ be a $V$-quantization of an operad $\P$. We say that $\P_V$
is coherent if $\P_V(n)$  is, for each
degree $n$, isomorphic  to $\P(n)\otimes V^{\otimes n}$
as a right $V^{\otimes n}$-module over the canonical isomorphism
$\epsilon(\P_V)\cong \P(n)$.
\end{definition}

Coherence in this sense measures the regularity of the behavior of the 
operadic ideal generated by $R_V$ in $\F_V(E)$.

Recall~(Definition~\ref{61}) 
that the collection of coherence constraints 
$\C$  of the operad $\P$ is defined as 
the indecomposables of the quotient 
$\ker (\pi)/\O$, where $\O$ is the
module of `obvious relations' of the operad $\P$. Let $p$
be the projection $\ker (\pi)\epic {\C}$. The main statement of
this section reads:

\begin{theorem}
\label{Tatouch}
The $V$-quantization $\P_V$ of $\P$ 
is coherent if and only if the composition
\begin{equation}
\label{zebrulka}
\Xi := \epsilon(\ker(\pi_V)) \stackrel{\rho}{\longrightarrow}
\ker(\pi)\stackrel{p}{\epic} {\C} 
\end{equation}
is an epimorphism.
\end{theorem}

\begin{proof}
The inclusion $\bfk \to V$, $\alpha \mapsto \alpha \cdot 1$,  
induces an embedding $\F(E) \hookrightarrow
\F_V(E)$ which is clearly a map of operads. 
This induces on $\F_V(E)\langle R
\rangle$, and thus also on $\ker(\pi_V)$, an $\F(E)$-module
structure. It is easy to conclude from Lemma~\ref{Nunyk} that
$\epsilon(\ker (\pi_V))$ is also a natural $\F(E)$-module 
and that the map $\rho: \epsilon(\ker (\pi_V) ) \to \ker(\pi)$ is
an $\F(E)$-homomorphism.

\begin{claim}
\label{Cudlik}
The map $\Xi : \eps(\ker(\pi_V)) \to \C$ 
of~(\ref{zebrulka}) is an epimorphism if and only if the
canonical map
$\rho : \epsilon(\ker (\pi_V) ) \to \ker(\pi)$ is onto.
\end{claim}

\noindent 
{\it Proof of the claim.\/}
If $\rho$ is onto, then clearly $\Xi$ must be an epimorphism, too.
To prove the opposite implication, observe that
a map $f: X \to Y$ of $\F(E)$-modules is an epimorphism if and only if
the composite of $f$ with the canonical 
projection $Y \epic Q_{\F(E)}(Y)$ to
the space of indecomposables is onto. 
This means that if $\Xi$ is an epimorphism,
then the composition $\epsilon(\ker(\pi_V)) \stackrel{\rho}{\longrightarrow}
\ker(\pi) \longrightarrow \ker(\pi)/\O$ is an epimorphism, too. 
As it is easy to verify, the map $\rho$ is always
an epimorphism on the space of
obvious relations, $\O \subset {\rm Im}{(\rho)}$, this means that the
map $\rho$ is onto. The claim is proved.

Summing up our observations, we have, for each $n\geq 1$, 
the sequence~(\ref{52q}) of right $\Vn$-modules
\[
0  \longrightarrow  \ker(\pi_V)(n)
\longrightarrow
\F_V(E)\langle R\rangle(n)
\stackrel{\pi_V}{\longrightarrow}
\F_V(E)(n)
\longrightarrow \P_V(n)  \longrightarrow 0
\]
and a sequence of $\bfk$-modules
\[
0
\longrightarrow \ker(\pi)(n)
\longrightarrow
\free{\F(E)}R(n)
\stackrel\pi\longrightarrow
\F(E)(n)
\longrightarrow \P(n)  \longrightarrow 0.
\]
We also know that $\eps(\F_V(E)\langle R\rangle(n))\cong
\free{\F(E)}R(n)$, $\eps(\F_V(E)(n))\cong \F(E)(n)$ 
and that $\eps\pi_V=
\pi\eps$. To finish the proof of Theorem~\ref{Tatouch}, 
it is, by Claim~\ref{Cudlik}, enough to
show that, for each $n\geq 1$,

\begin{equation}
\begin{minipage}{14cm}
\baselineskip18pt
$\P_V(n) \cong \P(n)\ot \Vn$ over the canonical isomorphism
$\eps(\P_V(n)) \cong \P(n)$ if and only if the canonical
map $\rho: \eps(\ker(\pi_V(n)))\to \ker(\pi(n))$ is an epimorphism. 
\end{minipage}
\end{equation}

\noindent 
This will clearly follow from the following lemma,
in which we put $W : = \Vn$.

\begin{lemma}
\label{zebrulinka}
Suppose we have an exact sequence
\begin{equation}
\label{klice}
0 \lra S \stackrel{\eta}{\lra} C \stackrel{\pi}{\lra} B 
\stackrel{\tau}{\lra} A \lra 0
\end{equation}
of $\bfk$-modules, and an exact sequence
\begin{equation}
\label{fotka}
0 \lra S' \stackrel{\eta'}{\lra} C' 
\stackrel{\pi'}{\lra} B' \stackrel{\tau'}{\lra} A' \lra 0
\end{equation}
of right $W$-modules. Suppose that the modules $C'$ and $B'$ are
$W$-free,
$\eps(C') \cong C$, $\eps(B')\cong
B$, and that, under these identifications, $\pi = \eps(\pi')$.

The canonical map $\rho : \eps(S') \to S$ is always a monomorphism and
the following two conditions are equivalent:
\begin{itemize}
\item[(i)]
the map $\rho : \eps(S') \to S$ is an epimorphism,
\item[(ii)]
$A' \cong A \ot W$ over the canonical isomorphism $\eps(A') \cong A$. 
\end{itemize}
\end{lemma}

\noindent 
{\it Proof of the lemma.\/}
To prove that the map $\rho$ is a monomorphism, observe that, since
$C'$ is a free $W$-module, the induced map $\eps(S') \to \eps(C')$ is
monic. Because $\rho$ is the composition of this map with the
identification $\eps(C')\cong C$, it must be a monomorphism as well.

For the rest of the proof we
may clearly assume there is a $\bfk$-vector space $D$ such
that~(\ref{klice}) is of the form
\[
0
\lra S \stackrel{\eta}{\lra} D \op S
\sta{\pi}{\lra} D\op A
\sta{\tau}{\lra} A \lra 0
\]
with $\eta(r) = (0,r)$, $\pi(d,r) = (d,0)$ and $\tau(d,a)=a$.
Then~(\ref{fotka}) is necessarily of the form
\[
0
\lra S' \stackrel{\eta'}{\lra} D \ot W\op S \ot W
\sta{\pi'}{\lra} D\ot W\op A \ot W
\sta{\tau'}{\lra} A' \lra 0
\]
with $\pi'$ represented by the matrix
\begin{equation}
\label{vifon}
\left(\begin{array}{cc}x&y\\z&w\end{array}\right),
\end{equation}
where
$x\in\HOM{D\ot W}{D\ot W}$, $y\in \HOM{S\ot W}{D\ot W}$,
$z\in \HOM{D\ot W}{ A\ot W}$ and
$w\in \HOM{S\ot W}{ A\ot W}$.
The assumption $\eps(\pi') = \pi$ translates to
\[ 
(\id_D\ot\epsilon)x=\id_D\ot \epsilon,\ 
(\id_D\ot\epsilon)y=0,\
(\id_A\ot\epsilon)z=0
\mbox { and }(\id_A\ot\epsilon)w=0
\]
and completeness of $W$ implies that the map $x$ is invertible.

In order to finish the proof, it is clearly enough to show that either
(i) or (ii) is equivalent to the existence of isomorphisms 
$\phi$ and $\psi$ of right $W$-modules 
such that the following diagram commutes:
\begin{equation}
\label{tool}
\begin{array}{ccc} 
D\ot W \oplus S\ot W 
&\stackrel{\pi\ot \id_W}{\verylra}&
D\ot W\oplus A\ot W
\\
\phi\uparrow&&\psi\uparrow
\\
D\ot W\oplus S\ot W
&\stackrel{\pi'}{\verylra}&
D\ot W\oplus A\ot W
\end{array}
\end{equation}
The existence of $\phi$ and $\psi$ and the commutativity of~(\ref{tool}) clearly
implies (i) and (ii) so it is enough  to prove the converse. Assuming either (i)
or (ii) we need to prove that there exist $\phi,\psi$ satisfying the equation
\begin{equation}
\label{Za}
\psi^{-1}\circ (\pi\ot\id_W)\circ\phi=\pi'.
\end{equation}
Let us suppose (i), that the map $\rho$ is an epimorphism. Put
\[
\phi:=\left(
\begin{array}{cc}x&y\\0&\id_{S\ot W}
\end{array}
\right)
\mbox{ and }
\psi:=\left(
\begin{array}{cc}\id_{D\ot W}&0\\-zx^{-1}&\id_{A\ot W}
\end{array}
\right).
\]
The map $\psi$ is clearly invertible and
\begin{eqnarray*}
\psi^{-1}\circ (\pi \ot \id)
\circ\phi
&=&
\left(
\begin{array}{cc}\id_{D\ot W}&0
\\
zx^{-1}&\id_{A\ot W}\end{array}
\right)
\left(
\begin{array}{cc}
\id_{D\ot W}&0\\0&0
\end{array}\right)
\left(
\begin{array}{cc}x&y\\0&\id_{S\ot W}
\end{array}
\right)
\\
&=&
\left(
\begin{array}{cc}x&y\\
z&zx^{-1}y\end{array}
\right).
\end{eqnarray*}
To prove that this composition gives $\pi'$ 
we need to show that $w=zx^{-1}y.$
Let $\{s_1,\ldots,s_n\}$ be a $\k$-basis for $S$. By the
surjectivity of $\rho$, we
know that for each $s_i$ there exists a corresponding $s'_i\in S'$,
such that $\rho(\eps(s'_i))= s_i$. 
Since $S'\subset D\ot W\oplus S\ot W$,
we can represent, for $1\leq i \leq n$, $s'_i$ as
$r'_i=d'_i +\sum_{1\leq j \leq n} s_j\ot w_{ij}$ with $w_{ij}\in W$ and
$d'_i\in D\ot W$. Then $\rho(\eps(s'_i))=s_i$
implies $\epsilon(w_{ij})=\delta_{ij}$, and the completeness of $W$
assures that the matrix $(w_{ij})$ is invertible. Thus $S'$ contains
elements of the form $s_i''=d_i''+ s_i\ot 1.$ For $s'\in S'$ there
exist $w_i$ such that $s'-\sum s''_i w_i\in S'\cap (D\otimes W).$
The invertibility of
the map $x$ implies that  $S'\cap (D\ot W)=0$, so the elements 
$\{s_i''\}_{1\leq i \leq n}$ span $S'$ as a $W$-module. 
They are clearly linearly independent over $W$ so they form  a basis.
The definition of $S'$ as the kernel of $\pi'$ gives the matrix equation
\begin{equation}
\label{mat}
\left(
\begin{array}{cc}x&y\\
z&w\end{array}
\right)
\left(
\begin{array}{c}u\\ \id_{S\ot W}
\end{array}
\right)
=
\left(
\begin{array}{c}0\\ 0
\end{array}
\right),
\end{equation}
where $u\in\HOM{S\ot W}{D\ot W}$ is defined by $u(s_i\ot 1)=d_i''$.
It follows from~(\ref{mat}) that
\[ xu+ y=0,\mbox{ and } zu+ w=0.\] 
Solving the first equation, $u=-x^{-1}y$, and substituting in the second
we obtain $w=zx^{-1}y$.

Suppose now~(ii). This means that there exists a 
$W$-isomorphism $\xi: A' \to A\ot W$
such that $(\id_A \ot \eps)\xi= \eps$. Consider the composition $\xi
\circ \tau' : D\ot W \op A \ot W \to A\ot W$ represented by the matrix
$(r,s)$, $r : D\ot W \to A\ot W$ and $s : A\ot W \to A \ot W$. By the
assumption on $\xi$, $(\id_A \ot \eps)s = \id_A \ot \eps$. So $s$ is
invertible and we can replace $\xi$ by $s^{-1} \circ \xi =: 
\overline\xi$. The
matrix representing the composition $\overline \xi \circ \tau$ 
is $(\overline r,\id_{\id \ot W})$. 

Now $\tau' \circ \pi' = 0$, so $\overline\xi\circ\tau '\circ \pi' = 0$ and
representing $\pi'$ as in~(\ref{vifon}) we conclude that
\[
(\overline r,\id_{A\ot W}) 
\left(
\begin{array}{cc}
x&y\\z&w
\end{array}
\right)
=
(0,0).
\]
Then
\[
\left(
\begin{array}{cc}
\id_{D\ot W}&0\\\overline r&\id_{A\ot W}
\end{array}
\right)
\left(
\begin{array}{cc}
x&y\\z&w
\end{array}
\right)
=
\left(
\begin{array}{cc}
x&y\\0&0
\end{array}
\right)
=
\left(
\begin{array}{cc}
\id_{D\ot W}&0\\0&0
\end{array}
\right)
\left(
\begin{array}{cc}
x&y\\0&\id_{R\ot W}
\end{array}
\right)
\]
and we may in~(\ref{tool}) take
\[
\phi :=\left(
\begin{array}{cc}
x&y\\0&\id_{R\ot W}
\end{array}
\right)
\mbox { and }
\psi:=\left(
\begin{array}{cc}
\id_{D\ot W}&0\\\overline r&\id_{A\ot W}
\end{array}
\right).
\]
The lemma, and thus also Theorem~\ref{Tatouch}, is proved.
\end{proof}

\begin{example}
\label{cervena_propiska}
{\rm
Let us investigate the coherence of the operad $\Ass_\Phi$ which we
introduced in Example~\ref{A}.
Clearly, $R_{{\Ass_\Phi}}$ is a free right $V^{\ot 3}$-module. 
Definition~(\ref{leftV}) of the left $V$-action gives
\begin{eqnarray*}
v \cdot \tilde r_{{\Ass_\Phi}} &=&
(1(23)) \cdot (\id \ot \Delta)\Delta(v) \cdot \Phi
- ((12)3) \cdot (\Delta \ot \id)\Delta(v)
\\
&=& \tilde r_{{\Ass_\Phi}} \cdot (\Delta \ot \id)\Delta(v) -
(1(23))[\Phi \cdot (\Delta \ot \id)\Delta(v) - (\id \ot \Delta)\Delta(v)\Phi].
\end{eqnarray*}
We see that $R_{{\Ass_\Phi}}$ is left $V$-closed if, for each $v\in V$,
\[
(\id \ot \Delta)\Delta(v)\cdot \Phi=\Phi\cdot(\Delta \ot \id)\Delta(v).
\]

We already know that the coherence constraints $\C$ of the
operad $\Ass = \eps(\Ass_\Phi)$ form an one dimensional 
$\bfk$-vector space
$\C = \C(4)$ with the generator corresponding to the
pentagon. Theorem~\ref{Tatouch} in this case says that the operad
$\Ass_\Phi$ is coherent if an only if the map $\rho :
\eps(\ker(\pi_{{\Ass_\Phi}})(4)) \to \ker(\pi_\Ass)(4)$ is onto. Since,
by Lemma~\ref{zebrulinka}, 
the map $\rho$ is a monomorphism, the later is
true if and only if  $\eps(\ker(\pi_{{\Ass_\Phi}})(4))\not =0$.

The `quantized' map $\pi_{\Ass_\Phi}:
\F_V(E)\langle R_{{\Ass_\Phi}} \rangle(4) \to (R_{{\Ass_\Phi}})(4)$
is described
by the following matrix with entries in $W= V^{\ot 4}$ (notation of
Example~\ref{hedges}): 
\begin{equation}
\label{oslicek}
\begin{tabular}{r|ccccc}
   &${\bf a}$&${\bf b}$&${\bf c}$&${\bf d}$&${\bf e}$\\ \hline
$\pi_{{\Ass_\Phi}}({\bf 1})$&$-1$&$(\Delta \OO \id^2)(\Phi)$
            &$0$&$0$&$0$\\
$\pi_{{\Ass_\Phi}}({\bf 2})$&$0$&$-1$&$
           (\id^2\OO \Delta)(\Phi)$&$0$&$0$\\
$\pi_{{\Ass_\Phi}}({\bf 3})$&$0$&$0$&$(1\OO \Phi)$&$
           -1$&$0$\\
$\pi_{{\Ass_\Phi}}({\bf 4})$&$0$&$0$&$0$&$(\id\OO\Delta\OO\id)\Phi$&$
           -1$\\
$\pi_{{\Ass_\Phi}}({\bf 5})$&$-1$&$0$&$0$&$0$&$(\Phi\OO 1)$
\end{tabular}
\end{equation}

Consider the kernel of this map. 
It follows from a very special form of the matrix for
$\pi_{\Ass_\Phi}(4)$
that $\pi_{\Ass_\Phi}({\bf 1}\cdot x_0 + {\bf 2}\cdot x_1
+\cdots + {\bf 5}\cdot x_4)=0$
can be expanded to the following system of equations for $\rada
{x_0}{x_4}\in W$:
\[
\alpha_i x_{\bar i}= x_{\bar i+1}\quad\mbox{for $\bar i\in
{\mathbb Z}_5$},
\]
with $\alpha_0=(\Delta \ot \id^2)(\Phi)$, $\alpha_1= (1\ot
\Phi)^{-1}(\id^2 \ot \Delta)(\Phi)$, $\alpha_2 = (\id \ot \Delta \ot
\id)(\Phi^{-1})$, $\alpha_3 = (\Phi \ot 1)^{-1}$ and $\alpha_4 =-1$.
We see that
\[
\ker(\pi_{\Ass_\Phi})(4) \cong
\{x \in W;\ (\alpha_4\cdots\alpha_0)x=x\}. 
\]

Clearly $\eps(\ker(\pi_{\Ass_\Phi})(4)) 
\not = 0$ if and only if
there
exists $x\in W$ such that $\eps(x) \not = 0$ and
$(\alpha_4\cdots\alpha_0)x=x$. By the completeness of $W$, 
such $x$ is invertible and
$(\alpha_4\cdots
\alpha_0)$ must equal $1$, which is the standard
pentagon identity for $\Phi$:
\begin{equation}
\label{PPP}
(1\ot\Phi)(\id \ot \Delta \ot\id)(\Phi)
(\Phi\ot 1)=(\id^2 \ot \Delta)(\Phi)(\Delta \ot \id^2)(\Phi). 
\end{equation}

We see that the coherence of the $V$-relative operad
$\Ass_\Phi$ means that the object
$V = (V,\cdot,\Delta,1,\Phi)$ 
is a quasi-bialgebra in the sense
of~\cite{drinfeld:Alg.iAnaliz89}. A natural example of an
$\Ass_\Phi$-algebra is the dot-construction 
of~\cite{markl-stasheff:JofAlg94}.%
\endex}\end{example}

\begin{example}
\label{modra_propiska}
{\rm\
Next consider the coherence of the generalized Lie operad.
Let $E$ and $\zeta$ have the same meaning as in Example~\ref{Lie} and
let $t\in \F(E)(3)$ be the element $\zeta(\zeta\ot \id)$,
corresponding to the bracketing $[[12]3]$ and $s$ the element $\zeta(\id \ot \zeta)$
corresponding to the bracketing $[1[23]]$. 
Define a symmetry condition on $\zeta$ by
 \begin{equation}\label{sym2}
\zeta\circ T_{21}\otimes \mathcalr =-\zeta, 
\end{equation}
where $$\mathcalr:=\sum \mathcalr_{1}\otimes \mathcalr_{2}\in V^{\ot 2}.$$
In other words, $E_V:=\F_V(E)(2)$ is the free rank one right $V^{\ot 2}$-module
which is the quotient of the rank two module with
 basis $\{\zeta, \zeta\circ T_{21}\}$ by the submodule with basis
$\{\zeta\circ T_{21}\otimes \mathcalr +\zeta\}.$
{}From identity~(\ref{leftV}) we have
$$
v\cdot(\zeta\circ T_{21}\otimes \mathcalr) =\zeta\circ T_{21}\otimes 
\Delta^{op}(v) \mathcalr $$
and
$$-v\cdot\zeta=-\zeta\ot\Delta(v)=\zeta\circ T_{21}\otimes \mathcalr \Delta(v).$$
Therefore the condition that $\F_V(E)(2)$ is a free right $V^{\ot 2}$-module
implies 
\begin{equation}\label{rconj}
\Delta^{op}(v)= \mathcalr \Delta(v)\mathcalr^{-1}.
\end{equation}
Iterating (\ref{sym2}), the free module condition implies Drinfel'd's
triangularity condition on $\mathcalr$:
\begin{equation}
\mathcalr_{21}\cdot \mathcalr=1.
\end{equation}
Operadic composition of equation~(\ref{sym2}) with
$\id\otimes \zeta$  implies the identity 
\begin{equation}\label{sym3}
t\circ T_{231}\ot (\id \ot \Delta)\mathcalr=-s.
\end{equation}
The ($\mathcalr,\Phi$)-relative version of the Lie operad is
defined by the basic tertiary relation
\begin{equation}\label{Lie-Phi}
\tilde r_{\Lie_{\mathcalr,\Phi}}:=t\ot 1 -s\otimes 
\Phi +s\circ T_{213}\otimes \Phi_{213}\mathcalr,
\end{equation}
which is our formalism for the `quantum Jacobi identity'
\begin{equation}
\label{qLie}
[[x,y],z]=\Phi[x,[y,z]]-\Phi\mathcalr_{21}[y,[x,z]].
\end{equation}

Let $R_{\Lie_{\mathcalr,\Phi}}$ be the right $V^{\ot 3}$-module generated by
$\tilde r_{\Lie_{\mathcalr,\Phi}}$. 
Requiring $\epsilon(\tilde r_{\Lie_{\mathcalr,\Phi}}) = r_{\Lie}$ (the
Jacobi identity) means that
\begin{equation}
\label{Kuratko}
\epsilon(\Phi) = 1
\mbox {\ and \ }
\epsilon(\mathcalr) =1.
\end{equation}
Again from the identity~(\ref{leftV}) we have
\begin{eqnarray*}
v\cdot\tilde r_{\Lie_{\mathcalr,\Phi}}&=&v\cdot(t\ot 1 -s\otimes 
\Phi +s\circ T_{213}\otimes \Phi_{213}\mathcalr)\\
&=&t\ot (\Delta\ot\id)\Delta(v)-s\ot 
(\id\ot \Delta)(\Delta(v))\Phi\\&&+s\circ T_{213}\ot 
T_{213}\left((\id\ot \Delta)\Delta(v)\right)\Phi_{213}\mathcalr.
\end{eqnarray*}
Together with the condition $R_V(3)$ be
a free right $V^{\ot 3}$-module this implies
\begin{eqnarray}\label{phicoassoc}
(\id\ot \Delta)(\Delta(v))\Phi&=&\Phi(\Delta\ot\id)\Delta(v)\mbox{ and }\\
T_{213}(\id\ot \Delta)(\Delta(v))\Phi\mathcalr_{21})&=&
\Phi_{213}\mathcalr(\Delta\ot\id)\ot\Delta(v).\label{phicoassoc2}
\end{eqnarray} 
Equation~(\ref{phicoassoc2}) follows from equation~(\ref{phicoassoc}) 
and equation~(\ref{rconj}).  

Next we investigate the symmetry condition~(\ref{sym3}). It is most
convenient to study these relations in the form of the quantum Jacobi
identity~(\ref{qLie}).
\begin{eqnarray}
[[x,y],z]&=&[\Phi_1 x,[\Phi_2 y,\Phi_3 z]]-
[\Phi_1 \mathcalr_2y,[\Phi_2\mathcalr_1 x,
\Phi_3 z]]\nonumber
\\
&=&-[\Phi_1 x,[\mathcalr_2\Phi_3 z,\mathcalr_1\Phi_2 y]]-
[\Phi_1 \mathcalr_2y,[\Phi_2\mathcalr_1 x,
\Phi_3 z]]\nonumber
\\
\label{Kacirek}
&=&-[\Phi_1 \mathcalr_2y,[\Phi_2\mathcalr_1 x,
\Phi_3 z]]-[[(\Phi^{-1})_1\Phi_1x,(\Phi^{-1})_2\mathcalr_2\Phi_3 z],
(\Phi^{-1})_3\mathcalr_1\Phi_2 y]\label{consist}\\
&&-
[\Phi_1\mathcalr_2(\Phi^{-1})_2\mathcalr_2\Phi_3z,
[\Phi_2\mathcalr_1(\Phi^{-1})_1\Phi_1x,
\Phi_3(\Phi^{-1})_3\mathcalr_1\Phi_2 y].\nonumber
\\
\nonumber 
&=&- \Phi\mathcalr_{21}[y,[x,z]] - \Phi^{-1}\mathcalr_{32}\Phi_{132}[[x,z],y]
- \Phi\mathcalr_{21}\Phi^{-1}_{213}  \mathcalr_{31}\Phi_{231} [z,[x,y]].
\end{eqnarray}
{}From the symmetry relation~(\ref{sym2}) we conclude:
\begin{eqnarray*}
\lefteqn{
(\Phi\mathcalr_{21}-(\id\otimes \Delta)\mathcalr_{21}\Phi_{231}^{-1}
\mathcalr_{13}\Phi_{213})[y,[ x,z]]-}
\\
&&-
((\id\otimes\Delta)\mathcalr_{21}
-\Phi\mathcalr_{21}\Phi^{-1}_{213}\mathcalr_{31}\Phi_{231})[z,[x,y]]=0.
\end{eqnarray*}
In order that $\left(\F(E)_V/R_V\right)(3)$ be a free $V^{\ot 3}$-module
of rank $2$ we must have the identity: 
\begin{equation}
\label{phihex}
(\De \ot \id)({\mathcal R}) =
\Phi_{312}{\mathcal R}_{13}\Phi^{-1}_{132}{\mathcal R}_{23}\Phi.
\end{equation}
This is the first Drinfel'd hexagon
identity, see \cite[3.9a, 3.9b]{drinfeld:Alg.iAnaliz89}.
The second hexagon identity 
\begin{equation}
\label{phiex2}
(\id \ot \De)({\mathcal R}) =
\Phi_{231}^{-1} {\mathcal R}_{13}\Phi_{213}{\mathcal R}_{12}\Phi^{-1}.
\end{equation}
 follows from the first one by the triangularity 
$\calR_{21} = \calR^{-1}_{12}$.

Next we want to determine the conditions which
guarantee the coherence of the operad $\Lie_{\mathcalr,\Phi} := 
\F_V(E)/(R_{\Lie_{\mathcalr,\Phi}})$
in the sense of Definition~\ref{32}.

We saw in Example~\ref{Lie} that  
the coherence constraints for $\Lie
= \eps(\Lie_{\mathcalr,\Phi})$ form the one-dimensional trivial representation of
$\Sigma_4$,  it is thus enough to
investigate when $\eps(\ker(\pi_{\Lie_{\mathcalr,\Phi}})(4))$
is nontrivial.

After determining the appropriate coefficients in
$V^{\otimes 4}$ given by extending the relation~(\ref{Lie-Phi})
to brackets with
four terms, we obtain a `$\Phi$-matrix' representing the map
\[
\pi_{\Lie(\mathcalr, \Phi)}:\F_V(E)\langle R_{\Lie(\mathcalr, \Phi)}
\rangle(4)\rightarrow(R_{\Lie,
\Phi}) \subset\F_V(E)(4).
\]
We will not write the full $\Phi$-matrix
here since it is to big to fit on a page (it is obtained by decorating
the entries of the matrix in Figure~\ref{matrix}); 
for our purposes it is
enough to observe that the upper left $5\times 5$ submatrix of this
$\Phi$-matrix is:
\begin{equation}
\label{opicak_Fuk}
\begin{tabular}{r|ccccc}
   &${\bf a}$&${\bf b}$&${\bf c}$&${\bf d}$&${\bf e}$\\ \hline
$\pi_{\Lie_{\mathcalr,\Phi}}({\bf 1})$&$1$&$
	   -(\Delta \otimes \id^2)(\Phi)$&$0$&$0$&$0$\\
$\pi_{\Lie_{\mathcalr,\Phi}}({\bf 2})$&$0$&$1$&$
	   -(\id^2\otimes \Delta)(\Phi)$&$0$&$0$\\
$\pi_{\Lie_{\mathcalr,\Phi}}({\bf 3})$&$0$&$0$&$-(1 \otimes \Phi)$&$
	   1$&$0$\\
$\pi_{\Lie_{\mathcalr,\Phi}}({\bf 4})$&$0$&$0$&$0$&
	     $-(\id \otimes \Delta \otimes \id)(\Phi)$&$
	   1$\\
$\pi_{\Lie_{\mathcalr,\Phi}}({\bf 5})$&$1$&$0$&$0$&$0$&$-(\Phi\otimes 1)$
\end{tabular}
\end{equation}
Observe that this matrix coincides, up to the sign reversal, with the
matrix~(\ref{oslicek}).

The next step is to describe $\ker(\pi_{\Lie_{\mathcalr,\Phi}})(4)$ as
a $V^{\otimes 4}$-module.  It is clear from the form of the matrix
that, as in Example~\ref{Lie},
the system of equations for
$\Ker{\pi_{\Lie_{\mathcalr,\Phi}}}(4)$ has the
form
\[
\alpha_j x_{i(j)}=x_{i'(j)},\ 1\leq j \leq 15,
\]
where $j$ is the index for a column, i.e.~an edge of the graph,
and $i(j),i'(j)$ are the two vertices adjacent to that edge,
i.e., the two rows with non-zero entries in that column.
The consistency conditions have the form
\[
\beta _i x_i=x_i
\]
where $\beta_i$ is a product of the $\alpha$'s going around a
closed path with initial and terminal vertex~$i$, which implies the
pentagon identity~(\ref{PPP}) as in the previous example. Thus we
must have both the pentagon and hexagon identities, and coherence of
the operad $\Lie_{\mathcalr,\Phi}$ implies Mac Lane coherence. The converse is
clear since the coefficients of a 
row in the matrix for $\pi_{\Lie_{\mathcalr,\Phi}}$ 
represent natural transformations between the three bracketings appearing
in a Jacobi relation, and Mac Lane coherence implies the uniqueness of the  
natural transformation connecting any two bracketings. 

More precisely, 
let $T_{{\bf xy}}$ be the element of $V^{\ot 4}$ representing the natural
transformation from a bracketing $\bf x$ to a 
bracketing $\bf y$. If we fix one
bracketing, such as ${\bf a}$ in the table in Figure~\ref{matrix},
multiplying each row of $\pi_{\Lie_{\mathcalr,\Phi}}$ by a suitable factor we get
a new matrix with the two non-zero entries in column ${\bf i}$ given
by $\pm T_{{\bf ai}}$.
Thus there is the same 
kind of dependency among the rows as in the classical case.}\endex
 \end{example}
\begin{theorem}
\label{zajicek_Cmuchalek}
The operad $\Lie_{\mathcalr,\Phi}$ is coherent if and only if 
$(V,\cdot,\Delta,\Phi,{\mathcal R})$ is a triangular quasi-Hopf 
algebra in the sense of~\cite{drinfeld:Alg.iAnaliz89}.
\end{theorem}

\begin{example}
\label{nozik}
{\rm\
Let us consider algebras consisting of a graded vector space $U_*$ and a
trilinear degree $-1$ product $\{-,-,-\} : U_*^{\ot 3}\to U_*$ 
satisfying, for all homogeneous
$a,b,c,d,e \in U_*$, 
\begin{equation}
\label{opulinka}
\{\{a,b,c\},d,e\}+ \sign{|a|} \{a,\{b,c,d\},e\} 
\sign{|a|+|b|}+ \{a,b,\{c,d,e\}\} = 0.
\end{equation}
The above means that $(U_*,\{-,-,-\})$ is an $A(\infty)$-%
algebra~\cite{stasheff:TAMS63} 
with all structure operations $\mu_n$ trivial except $\mu_3 =
\{-,-,-\}$. Let ${\mathcal A}$ be the non-$\Sigma$ operad describing these
algebras. It is not quadratic, but it is homogeneous in the
sense that ${\mathcal A} = \prez ER$ with $E = E(3)$, the one dimensional
space generated
by the product $\{-,-,-\}$, and $R \subset
\F(E)(5)$ generated by the left hand side of~(\ref{opulinka}). 
For these homogeneous operads, it still
makes sense to introduce their $!$-duals, and an  easy
calculation gives
\[
{\mathcal A}^!(n) =
\cases {\downarrow^n \bfk}{for $n = 1$ (mod $2$), while}0{otherwise,}
\]
where $\downarrow^n \bfk$ is the one-dimensional graded vector
space concentrated in degree $-n$.
It can also be shown that ${\mathcal A}$ 
is Koszul in a suitably generalized sense,
thus $\C_{{\mathcal A}} = \C_{{\mathcal A}}(7)= {\mathcal A}^!(7)$ 
is one-dimensional. Therefore the
kernel of the map
\[
\pi(7) : \F(E)\langle R \rangle (7) \lra \F(E)(7)
\]
is also of dimension one. In fact, $\dim \F(E)\langle R \rangle (7)=8$,
$\dim \F(E)(7) = 12$ and $\pi(7)$ is represented by the matrix
 \begin{equation}
\label{kozulinka}
\left(
\begin{tabular}{rrrrrrrrrrrr}
$+1$&$+1$&$+1$&$ 0$&$ 0$&$ 0$&$ 0$&$ 0$&$ 0$&$ 0$&$ 0$&$ 0$
\\
$ 0$&$ 0$&$ 0$&$+1$&$+1$&$+1$&$ 0$&$ 0$&$ 0$&$ 0$&$ 0$&$ 0$
\\
$ 0$&$ 0$&$ 0$&$ 0$&$ 0$&$ 0$&$+1$&$+1$&$+1$&$ 0$&$ 0$&$ 0$
\\
$ 0$&$+1$&$ 0$&$ 0$&$ 0$&$ 0$&$ 0$&$ 0$&$ 0$&$-1$&$-1$&$ 0$
\\
$ 0$&$ 0$&$+1$&$ 0$&$ 0$&$+1$&$ 0$&$ 0$&$ 0$&$ 0$&$ 0$&$-1$
\\
$+1$&$ 0$&$ 0$&$+1$&$ 0$&$ 0$&$+1$&$ 0$&$ 0$&$ 0$&$ 0$&$ 0$
\\
$ 0$&$ 0$&$ 0$&$ 0$&$+1$&$ 0$&$ 0$&$+1$&$ 0$&$+1$&$ 0$&$ 0$
\\
$ 0$&$ 0$&$ 0$&$ 0$&$ 0$&$ 0$&$ 0$&$ 0$&$+1$&$ 0$&$+1$&$+1$
\end{tabular}
\right)
\end{equation}

\noindent 
Since~(\ref{kozulinka}) has exactly two nontrivial entries in
each column, the corresponding TA-structure is graphlike, with
Tel-A-graph the `M\"obius strip:' 
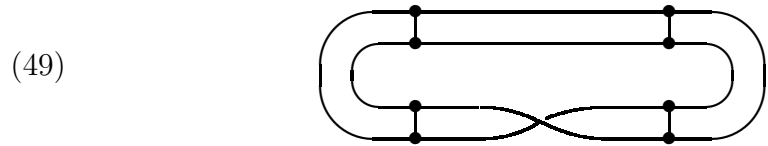
\begin{equation}
\label{mobius}
\def\qbezier{\bezier{100}}
{
\thicklines
\unitlength=1.200000pt
\begin{picture}(170.00,20.00)(0.00,20.00)
\put(110.00,0.00){\makebox(0.00,0.00){$\bullet$}}
\put(110.00,10.00){\makebox(0.00,0.00){$\bullet$}}
\put(30.00,0.00){\makebox(0.00,0.00){$\bullet$}}
\put(30.00,10.00){\makebox(0.00,0.00){$\bullet$}}
\put(110.00,30.00){\makebox(0.00,0.00){$\bullet$}}
\put(110.00,40.00){\makebox(0.00,0.00){$\bullet$}}
\put(30.00,30.00){\makebox(0.00,0.00){$\bullet$}}
\put(30.00,40.00){\makebox(0.00,0.00){$\bullet$}}
\put(110.00,40.00){\line(0,-1){10.00}}
\put(110.00,0.00){\line(0,1){0.00}}
\put(110.00,10.00){\line(0,-1){10.00}}
\put(30.00,10.00){\line(0,-1){10.00}}
\put(30.00,40.00){\line(0,-1){10.00}}
\qbezier(71.00,6.50)(79.75,10.00)(90.25,10.00)
\qbezier(50.25,0.00)(59.75,0.00)(68.00,4.75)
\qbezier(70.00,5.00)(80.00,0.00)(90.00,0.00)
\qbezier(50.50,10.00)(60.00,10.00)(69.75,5.00)
\put(90.00,0.00){\line(1,0){30.00}}
\put(90.00,10.00){\line(1,0){30.00}}
\put(20.00,0.00){\line(1,0){30.00}}
\put(20.00,10.00){\line(1,0){30.00}}
\put(20.00,30.00){\line(1,0){100.00}}
\put(20.00,40.00){\line(1,0){100.00}}
\put(120.00,20.00){\oval(40.00,40.00)[r]}
\put(20.00,20.00){\oval(40.00,40.00)[l]}
\put(120.00,20.00){\oval(20.00,20.00)[r]}
\put(20.00,20.00){\oval(20.00,20.00)[l]}
\end{picture}}
\end{equation}

\vskip 7mm
\noindent 
with $8$ vertices corresponding to the rows of~(\ref{kozulinka}) and
12 edges corresponding to the columns of~(\ref{kozulinka}). 

Let $V
= (V,\cdot,1,\eps)$ be a unital augmented algebra
and let $\Delta = \Delta^3 : V \to
V^{\ot 3}$ be a homomorphism. Let us `quantize'~(\ref{opulinka}) by
decorating it by some invertible $\Phi, \Psi \in V^{\ot 3}$ to
\[
\{\{a,b,c\},d,e\}+ \sign{|a|}\Phi \{a,\{b,c,d\},e\}+ \sign{|a|+|b|} 
\Psi\{a,b,\{c,d,e\}\} = 0.
\]
It is an easy exercise to prove that the left $V$-invariance of this
axiom implies that
\[
(\id \ot \Delta \ot \id)\Delta \cdot \Phi = \Phi\cdot 
(\Delta \ot \id^2)\De\
\mbox { and }
(\id^2 \ot \Delta)\Delta \cdot \Psi = \Psi\cdot (\Delta \ot \id^2)\Delta.
\]
The coherence of the corresponding $V$-relative operad 
${\mathcal A}_{\Phi,\Psi}$ then means that
$\Phi$ and $\Psi$ satisfy 6 equations corresponding to 6 generators of
the fundamental group of~(\ref{mobius}). 
To get these equations would mean to evaluate the `decorated' matrix for
$\pi_{\Phi, \Psi}(7)$ from which these equations can be easily read
off. We leave this to the reader.
\endex}
\end{example}


%

\end{document}